\documentclass[11pt]{article}
\usepackage{amssymb,amsmath}

\newtheorem{theorem}{Theorem}[section]
\newtheorem{lemma}[theorem]{Lemma}
\newtheorem{proposition}[theorem]{Proposition}

\newtheorem{definition}[theorem]{Definition}
\newtheorem{corollary}[theorem]{Corollary}
\newtheorem{remark}[theorem]{Remark}

\begin{document}

\title{Isomorphic Schauder decompositions in certain Banach spaces
}

\author{ Vitalii Marchenko\thanks{V.N. Karazin Kharkiv National University, Department of Mathematics and\newline Mechanics, Svobody Sq. 4, 61022, Kharkiv, Ukraine, \tt vitalii.marchenko@karazin.ua}
}

\maketitle

\begin{abstract}
  We extend a theorem of Kato on similarity for sequences of projections in Hilbert spaces to the case of isomorphic Schauder decompositions in certain Banach spaces.
To this end we use $\ell_{\Psi}$-Hilbertian and $\infty$-Hilbertian Schauder decompositions instead of orthogonal Schauder decompositions, generalize the concept of an orthogonal Schauder decomposition in a Hilbert space and introduce the class of spaces with Schauder-Orlicz decompositions.
Furthermore, we generalize the notions of type, cotype, infratype and $M$-cotype of a Banach space and study the properties of unconditional Schauder decompositions in spaces possessing certain geometric structure.
\end{abstract}

\section{Introduction}

Throughout the paper $E$ will denote a Banach space over the field $\mathbb{R}$ or $\mathbb{C}$ and $\mathbb{Z}_{+}$ will denote a set of nonnegative integers. A sequence $\{\mathfrak{M}_n\}_{n=0}^{\infty}$ of (not necessarily closed) nonzero
linear subspaces of $E$ is called a decomposition (or basis of subspaces) of $E$ provided every element $x\in E$ has a unique, norm convergent
expansion $x=\sum\limits_{n=0}^{\infty} x_n,$ where $x_n\in \mathfrak{M}_n$ for $n=0,1,\dots$ In other words, $\{\mathfrak{M}_n\}_{n=0}^{\infty}$ is
a decomposition of $E$ if and only if $E$ is an infinite direct sum of subspaces $\{\mathfrak{M}_n\}_{n=0}^{\infty}$:
$E=\bigoplus\limits_{n=0}^{\infty} \mathfrak{M}_n$ \cite{Grinblyum}. If $\{\mathfrak{M}_n\}_{n=0}^{\infty}$ is a decomposition of $E$, the sequence
of linear projections $\{P_n\}_{n=0}^{\infty}$ on $E$, defined by $P_n x=x_n,$ $n=0,1,\dots$, where $x=\sum\limits_{n=0}^{\infty} x_n,$ $x_n\in \mathfrak{M}_n$, is called the sequence of coordinate projections associated to the decomposition $\{\mathfrak{M}_n\}_{n=0}^{\infty}$, or, shortly,
the associated sequence of coordinate projections (a.s.c.p.). Thus, if $\{\mathfrak{M}_n\}_{n=0}^{\infty}$ is a decomposition of $E$, then each $x\in E$ has a unique expansion $x=\sum\limits_{n=0}^{\infty} P_n x$, where $\{P_n\}_{n=0}^{\infty}$ is the a.s.c.p. The uniqueness condition of the definition of decomposition yields $P_i P_j=\delta_{i}^{j} P_i= \delta_{i}^{j} P_j$ for $i,j=0,1,\dots,$ that is, every a.s.c.p. consists of mutually
orthogonal projections.

If each $\mathfrak{M}_n$ is closed or, equivalently, if each of the
projections $P_n$ of the a.s.c.p. is continuous on $E$,
the decomposition is called a Schauder decomposition \cite{Singer2}. In case when every $\mathfrak{M}_n$ is finite dimensional, the decomposition is called a finite
dimensional decomposition, or simply an FDD. A Schauder basis can be considered as an FDD $\{\mathfrak{M}_n\}_{n=0}^{\infty}$ such that every subspace $\mathfrak{M}_n$ is one dimensional.
An FDD is at the same time a Schauder decomposition and an example of decomposition, which is
not a Schauder decomposition, may be found in \cite{Sanders2}.

The concept of Schauder decomposition is
a natural generalization of the Schauder basis concept and was first introduced by Grinblyum \cite{Grinblyum}. Independently, Fage \cite{Fage1,Fage2}
has studied this concept in Hilbert spaces. In the 1960's
Schauder decompositions and their various properties were studied by Marcus \cite{Marcus} and Kato \cite{Kato2} (Hilbert space case), Retherford
\cite{Retherford1,Retherford3}, Sanders \cite{Sanders1,Sanders2}, Davis \cite{Davis} and others. The concept of $\pi_{\lambda}$ and dual $\pi_{\lambda}$ space was introduced by Johnson \cite{Johnson1} and its direct relevance to the existence of FDD's was found. More precisely, a separable Banach space $E$ has an FDD if and only if $E$ is a dual $\pi_{\lambda}$ space.
A Schauder decomposition $\{\mathfrak{M}_n\}_{n=0}^{\infty}$ is called unconditional if the expansions $x=\sum\limits_{n=0}^{\infty} x_n,$ $x_n\in \mathfrak{M}_n$, converge unconditionally for each $x\in E$. The Schauder decomposition $\{\mathfrak{M}_n\}_{n=0}^{\infty}$ is said to be conditional if it is non-unconditional. The general problem of characterizing conditions under which one can construct an unconditional basis for $E$ by forming an unconditional basis for each $\mathfrak{M}_n$, where $\{\mathfrak{M}_n\}_{n=0}^{\infty}$ is an unconditional FDD, was considered in \cite{Casazza}. One sufficient condition for the product of two a.s.c.p.'s, which correspond to unconditional Schauder decompositions, to be also the a.s.c.p. corresponding to unconditional Schauder decomposition was established in \cite{DePagter}.

Schauder decompositions together with Schauder bases are effective tools of functional analysis and techniques associated with them are widely applied in infinite dimensional linear systems theory, see, e.g., \cite{Curtain,Rabah1,Rabah2,Rabah3,Zwart}. An interrelation between Schauder decompositions, unconditional Schauder decompositions and the problems of the best approximation in a Banach space has been noted in \cite{Jain2,Jain3}.
An interaction between Schauder decompositions, the property of $R$-boundedness of collections of operators and geometric properties of a Banach space was studied in \cite{Clement}, and some relations between the Grothendieck and Dunford-Pettis properties in K\"{o}the echelon spaces of infinite order and Schauder decompositions have been noted in \cite{Bonet}. The discussion of an interplay between frequent hypercyclicity, chaos and unconditional Schauder decompositions is presented in \cite{DelaRosa}.

Although bases and decompositions have a lot of similar properties, there are results concerning decompositions which do not have analogues
for bases or whose analogues for bases are not known or are false.
E.g., every Banach space (actually, every normed linear space) has a decomposition \cite{Sanders2}, but $\ell_{\infty}$ has no Schauder decomposition \cite{Singer2},
and there even exists a separable Banach space with no Schauder decomposition \cite{Allexandrov}, whereas the classical basis question has only a negative answer (Enflo's counterexample).
Moreover, all spaces with bases should be separable while there are nonseparable Banach spaces which possess Schauder decompositions  \cite{Sanders2,Chadwick}. For more about Schauder decompositions and a.s.c.p.'s see, e.g., \cite{Bilalov2,Gohberg,Johnson2,Lindenstrauss,Singer2}.

The study of Schauder decompositions in Hilbert spaces was motivated by the problems raised by spectral theory, in particular, by the problem of unconditional basicity of the eigenprojections of a nonselfadjoint and unbounded, in general, linear operator. Marcus \cite{Marcus} uses his results to obtain conditions providing an unconditional basicity of eigen and root vectors (root subspaces) of a dissipative operator, see \cite{Gohberg} for details. In the 1967 Kato has presented one sufficient condition for the similarity of sequences of projections in Hilbert spaces.
\begin{theorem}[Kato \cite{Kato2}]\label{Kato}
Let $\{P_n\}_{n=0}^{\infty}$ be a sequence of nonzero selfadjoint projections in Hilbert space $H$ such that $\sum\limits_{n=0}^{\infty}P_n=I$, and $\{J_n\}_{n=0}^{\infty}$ be a sequence of nonzero projections in $H$ satisfying $J_n J_m=\delta_{n}^{m} J_n$ for $n,m=0,1,\dots$ Furthermore, assume that
\begin{equation}\label{1}
    \dim P_0 = \dim J_0 = m <\infty,
\end{equation}
\begin{equation}\label{2}
   \sum\limits_{n=1}^{\infty} \|P_n (J_n-P_n) x\|^2 \leq c^2 \|x\|^2 \quad \text{for all $x\in H$},
\end{equation}
where $c$ is a constant such that $0\leq c<1$. Then $\{J_n\}_{n=0}^{\infty}$ is similar to $\{P_n\}_{n=0}^{\infty}$, that is, there exists an isomorphism $S$ such that $J_n=S P_n S^{-1}$ for $n=0,1,\dots$
\end{theorem}
It follows from Theorem \ref{Kato} that the sequence of subspaces $\{\mathfrak{M}_n=P_n H\}_{n=0}^{\infty}$ generate an orthogonal Schauder decomposition of $H$, and the sequence $\{\mathfrak{N}_n=J_n H\}_{n=0}^{\infty}$ is also the Schauder decomposition of $H$ which is isomorphic to $\{\mathfrak{M}_n\}_{n=0}^{\infty}$.
This result gave a new impetus to further investigations in spectral theory. In the 1968 Clark \cite{Clark} applied Theorem \ref{Kato} to the study of spectral properties of relatively bounded perturbations of ordinary differential operators. In the 1972 Hughes \cite{Hughes} applied Theorem \ref{Kato} to prove some perturbation theorems for relative spectral problems.
Kato \cite{Kato1} considered the problem of completeness of the eigenprojections of slightly nonselfadjoint operator as a perturbation problem for a selfadjoint operator and based the solution of this problem on his Theorem \ref{Kato}. Recently, Adduci and Mityagin applied Theorem \ref{Kato} to the study of eigenfunction expansions of the perturbed harmonic oscillator $L=-\frac{d^2}{dx^2} + x^2 +B$, $B=b(x)$ with domain in $L_2(\mathbb{R})$ \cite{Adduci1}, and to analysis of the perturbation $A=T+B$ of a selfadjoint operator $T$ in Hilbert space $H$ with discrete spectrum \cite{Adduci2}. It was shown in \cite{Adduci1} that the eigensystem of $L$ forms an unconditional basis for $L_2(\mathbb{R})$ in case when $b$ belongs to certain function spaces (e.g., $L_p(\mathbb{R})$ with $p\in [2,\infty)$), and it was proved in \cite{Adduci2} that, under certain conditions on $T$ and $B$, eigen and root system of the operator $A$ constitutes an unconditional basis for $H$.

The main purpose of this paper is to extend the result of Kato (Theorem \ref{Kato}) for the case of certain sequences of projections $\{P_n\}_{n=0}^{\infty}$ in Orlicz sequence spaces $\ell_{\Phi}$, and for the case of certain sequences of projections $\{P_n\}_{n=0}^{\infty}$ in spaces with Schauder-Orlicz decompositions, which are introduced in Section 5. The second purpose is to study the properties of unconditional Schauder decompositions in spaces possessing certain geometric characteristics.

The paper is organized as follows.
Section 2 deals with the notion of isomorphic Schauder decompositions and gives a short survey. Section 3 provides intrinsic characterizations of unconditional Schauder decompositions in Banach and Hilbert spaces. In Section 4 the generalized concepts of type, cotype, infratype and $M$-cotype of a Banach space are introduced and an interplay between these concepts, $\ell_{\Phi}$-Hilbertian, $\ell_{\Psi}$-Besselian Schauder decompositions and unconditional Schauder decompositions is studied. Section 5 introduces the concept of Schauder-Orlicz decomposition and focuses on the main results of the paper concerning isomorphic Schauder decompositions in Orlicz spaces $\ell_{\Phi}$ and isomorphic Schauder decompositions in spaces with Schauder-Orlicz decompositions. Proofs are given in Section 6 and conclusions are drawn in Section 7.

\section{Isomorphic Schauder decompositions}
We say that Schauder decompositions $\{\mathfrak{M}_n\}_{n=0}^{\infty}$ and $\{\mathfrak{N}_n\}_{n=0}^{\infty}$ of $E$ are isomorphic provided there exists an isomorphism $S$ on $E$ such that 
\begin{equation}\label{1789876}
 \mathfrak{N}_n=S\mathfrak{M}_n,\: n=0,1,\dots  
\end{equation}
Note that any two sequences of subspaces (not necessary Schauder decompositions) $\{\mathfrak{M}_n\}_{n=0}^{\infty}$, $\{\mathfrak{N}_n\}_{n=0}^{\infty}$ of $E$ satisfying (\ref{1789876}) are called in \cite{Singer2} fully equivalent. Despite this we find our terminology more convenient.

Isomorphic Schauder decompositions was first considered in \cite{Fage1,Fage2} for Hilbert spaces. In fact, the paper \cite{Fage1} deals with a.s.c.p.'s, which are similar to associated systems of coefficient orthoprojections (a.s.c.o.), but, in view of Proposition \ref{proposition} below, it is the same. The main result of \cite{Fage1} is the criterion for the a.s.c.p. to be similar to the a.s.c.o. The following proposition asserts that two Schauder decompositions in a Banach space $E$ are isomorphic if and only if their a.s.c.p.'s are similar.
\begin{proposition}\label{proposition}
For Schauder decompositions $\{\mathfrak{M}_n\}_{n=0}^{\infty}$ and $\{\mathfrak{N}_n\}_{n=0}^{\infty}$ in $E$ with a.s.c.p.'s $\{P_n\}_{n=0}^{\infty}$ and $\{J_n\}_{n=0}^{\infty}$, respectively, the following statements are equivalent:
$$(i)\:\: \{\mathfrak{M}_n\}_{n=0}^{\infty}\:\: \text{and}\:\: \{\mathfrak{N}_n\}_{n=0}^{\infty} \:\: \text{are isomorphic};$$
$$(ii)\:\: \{P_n\}_{n=0}^{\infty}\:\: \text{and}\:\: \{J_n\}_{n=0}^{\infty} \:\: \text{are similar}.$$
\end{proposition}
In the 1960's, Marcus \cite{Marcus} and Kato \cite{Kato2} studied isomorphic Schauder decompositions in Hilbert spaces in connections with the problems of spectral theory. The study of isomorphic Schauder decompositions in Banach spaces is much more complicated. It is caused by the complexity and diversity of the geometry of Banach spaces. However, there are some interesting results in this direction. An equivalence of sequences of subspaces in Banach spaces and in dual spaces was studied in \cite{Jain1,Ahmad}. Certain results concerning basic sequences of subspaces, unconditional basic sequences of subspaces and the Paley-Wiener stability criterion were obtained in \cite{Retherford1}. Some other results concerning isomorphic Schauder decompositions may be found in \cite{Bilalov2,Singer2,Gohberg}.

\section{Unconditional Schauder decompositions}

The sufficient condition for the existence of (unconditional) Schauder decomposition in $E$ is the existence of complemented subspace $\mathfrak{L}\subset E$ with an (unconditional) basis \cite{Singer2}. Thus, e.g., $C[0,1]$ and $L_1[0,1]$ have unconditional Schauder decompositions, although they have no unconditional basis \cite{Singer1}.
Another remarkable distinction between spaces with unconditional basis and spaces possessing unconditional Schauder decompositions is the following.
There exist Banach spaces (e.g., $C[0,1]$ and $L_1[0,1]$) which cannot be embedded isomorphically into any Banach space with an unconditional basis \cite{Singer1}. Moreover, these spaces cannot be embedded into any Banach space with an FDD. However, for unconditional Schauder decompositions we have that every Banach space $E$ can be embedded, as a complemented subspace, into a Banach space with an unconditional Schauder decomposition, e.g., into $E\times \ell_1$ or $E\times c_0$ \cite{Singer2}.

All unconditional Schauder decompositions in a Banach space $E$ have the following intrinsic characterizations.
\begin{theorem}[\cite{Singer2,DePagter,Jain1}]\label{USD}
Let $\{\mathfrak{M}_n\}_{n=0}^{\infty}$ be a sequence of nonzero subspaces of $E$ such that $\overline{Lin} \{\mathfrak{M}_n\}_{n=0}^{\infty} = E$. Then the following statements are equivalent:

(1) $\{\mathfrak{M}_n\}_{n=0}^{\infty}$ is an unconditional Schauder decomposition of $E$.

(2) Every permutation $\{\mathfrak{M}_{\sigma(n)}\}_{n=0}^{\infty}$ of $\{\mathfrak{M}_n\}_{n=0}^{\infty}$ is a Schauder decomposition of $E$.

(3) For every increasing sequence $\{s_n\}_{n=0}^{\infty}\subset \mathbb{Z}_{+}$, the subspaces\newline $\overline{Lin} \{\mathfrak{M}_{s_n}\}_{n=0}^{\infty}$ and $\overline{Lin} \{\mathfrak{M}_{z_n}\}_{n=0}^{\infty}$, where $\{z_n\}_{n=0}^{\infty}$ is an increasing sequence of $\mathbb{Z}_{+}$ such that $\{s_n\}_{n=0}^{\infty} \cup \{z_n\}_{n=0}^{\infty} =\mathbb{Z}_{+}$, are complementary to each other, i.e.
$$\overline{Lin}
\{\mathfrak{M}_{s_n}\}_{n=0}^{\infty}\oplus \overline{Lin}  \{\mathfrak{M}_{z_n}\}_{n=0}^{\infty} = E.$$
(4) There exists a constant $M\geq 1$, such that
\begin{equation}\label{M}
    \left\|\sum\limits_{i=0}^{n} \delta_i y_i\right\| \leq M \left\|\sum\limits_{i=0}^{n}  y_i\right\|
\end{equation}
for any $n\in\mathbb{Z}_{+}$, $y_i\in \mathfrak{M}_i$ and $\delta_i\in \{0,1\}$ $(i=0,1,\dots,n).$

(5) There exists a constant $M_1\geq 1$, such that
$$\left\|\sum\limits_{i=0}^{n} \varepsilon_i y_i\right\| \leq M_1 \left\|\sum\limits_{i=0}^{n}  y_i\right\|$$
for any $n\in\mathbb{Z}_{+}$, $y_i\in \mathfrak{M}_i$ and $\varepsilon_i\in \{-1,1\}$ $(i=0,1,\dots,n).$

(6) There exists a constant $M_2\geq 1$, such that
$$\left\|\sum\limits_{i=0}^{n} \beta_i y_i\right\| \leq M_2 \left\|\sum\limits_{i=0}^{n}  y_i\right\|$$
for any $n\in\mathbb{Z}_{+}$, $y_i\in \mathfrak{M}_i$ and scalars $\beta_i$ with $|\beta_i|\leq 1$ $(i=0,1,\dots,n).$

(7) There exists a constant $M_3\geq 1$, such that
$$\left\|\sum\limits_{j=0}^{n} y_{i_j}\right\| \leq M_3 \left\|\sum\limits_{i=0}^{n}  y_{i_j} + \sum\limits_{j=0}^{m}  y_{k_j}\right\|$$
for any $n,m\in\mathbb{Z}_{+}$ and any $y_{i_j}\in \mathfrak{M}_{i_j}$ $(j=0,1,\dots,n)$, $y_{k_j}\in \mathfrak{M}_{k_j}$ $(j=0,1,\dots,m)$ which satisfy $\{i_0,\dots, i_n\} \cap \{k_0,\dots, k_m\}=\emptyset.$
\end{theorem}
\begin{definition}\label{M-USD}
We will say that Schauder decomposition $\{\mathfrak{M}_n\}_{n=0}^{\infty}$ is unconditional with constant $M$ if and only if condition \textit{(4)} of Theorem \ref{USD} holds.
\end{definition}
One can also define orthogonal and hyperorthogonal Schauder decompositions, by properties \textit{(7)} and \textit{(6)} of Theorem \ref{USD} with $M_3=1$ and $M_2=1$, respectively. It must be emphasized that, apart from the above, all unconditional Schauder decompositions in a Hilbert space $H$ have an additional intrinsic characterizations. We observe that a Schauder decomposition $\{\mathfrak{M}_n\}_{n=0}^{\infty}$ in $H$ is orthogonal if and only if all the subspaces $\mathfrak{M}_n$, $n=0,1,\dots$ are mutually orthogonal. Clearly, the a.s.c.p. for an orthogonal Schauder decomposition $\{\mathfrak{M}_n\}_{n=0}^{\infty}$ in $H$ is the a.s.c.o. Combining Proposition \ref{proposition} and Theorem \ref{USD} with some results of \cite{Fage1} and \cite{Wermer} we have the following.

\begin{theorem}\label{Orthogonal}
For Schauder decomposition $\{\mathfrak{N}_n\}_{n=0}^{\infty}$ in $H$ with the a.s.c.p. $\{J_n\}_{n=0}^{\infty}$ the following statements are equivalent:

(1) $\{\mathfrak{N}_n\}_{n=0}^{\infty}$ is an unconditional Schauder decomposition.

(2) There exists an orthogonal Schauder decomposition $\{\mathfrak{M}_n\}_{n=0}^{\infty}$ such that $\{\mathfrak{N}_n\}_{n=0}^{\infty}$ and $\{\mathfrak{M}_n\}_{n=0}^{\infty}$ are isomorphic.

(3) There exists the a.s.c.o. $\{P_n\}_{n=0}^{\infty}$ such that $\{J_n\}_{n=0}^{\infty}$ and $\{P_n\}_{n=0}^{\infty}$ are similar.

(4) There exists a constant $C\geq 1$ such that
\begin{equation}\label{2}
    \frac{1}{C} \|x\|^2 \leq \sum\limits_{n=0}^{\infty} \|J_n x\|^2 \leq C\|x\|^2 \quad \text{for every}\:\: x\in H.
\end{equation}
\end{theorem}
Observe that every sequence of projections $\{J_n\}_{n=0}^{\infty}$ satisfying (\ref{2}) is called in \cite{Zwart} a Riesz family.
Theorem \ref{Orthogonal} provides the uniqueness, up to an isomorphism, of unconditional Schauder decomposition in a Hilbert space $H$. And, consequently, one obtains the uniqueness of (bounded) unconditional basis in each separable Hilbert space $H$. On the other hand, a separable Hilbert space has uncountably many mutually non-isomorphic conditional bases \cite{Johnson3}. In fact, we have the following more general result. Namely, every Banach space possessing an infinite (Schauder) basis has uncountably many mutually non-isomorphic bases \cite{Johnson3}. In contrast to Hilbert spaces, each of the spaces $\ell_p$ and $L_p[0,1]$ for $p\in(1,\infty)\setminus {2}$ have two bounded unconditional bases which are not isomorphic \cite{Singer1}. Moreover, a Banach space has, up to an isomorphism, a unique unconditional basis if and only if it is isomorphic to one of the spaces $\ell_1$, $\ell_2$ or $c_0$ \cite{Johnson3}.

\section{Banach space geometry and isomorphic unconditional Schauder decompositions}

\subsection{Orlicz-Rademacher structural properties of Banach\newline spaces}
\begin{definition}[\cite{Lindenstrauss}]\label{Orlicz function}
An Orlicz function $\Phi$ is a continuous non-decreasing and convex function defined for $t\geq 0$ such that $\Phi(0)=0$ and $\lim\limits_{t \rightarrow \infty} \Phi(t)=\infty$.
\end{definition}
To any Orlicz function $\Phi$ one can associate the space $\ell_{\Phi}$ of all sequences of (complex) scalars $x=(a_0,a_1,\dots)$ such that $\sum\limits_{n=0}^{\infty} \Phi\left(\frac{|a_n|}{\rho}\right) < \infty$ for some $\rho>0$. The space $\ell_{\Phi}$ equipped with the norm
$$\|x\|_{\ell_{\Phi}}=\inf\left\{\rho >0:\quad  \sum\limits_{n=0}^{\infty} \Phi\left(\frac{|a_n|}{\rho}\right)\leq 1\right\}$$
is a Banach space usually called an Orlicz sequence space. One can see that, for $\Phi(t)=t^p$ with $p\geq 1$, $\ell_{\Phi}$ coincides with the classical space $\ell_p$. Consider $\{x_j\}_{j=0}^{n}\subset E$, $n\in\mathbb{N}$. In what follows we will use the notation
$$\mathbb{E}\left\|\sum\limits_{j=0}^{n} \varepsilon_j x_j \right\| =
\int\limits_{0}^{1} \left\|\sum\limits_{j=0}^{n} r_j (t) x_j
\right\| dt = \frac{1}{2^n}\sum\limits_{\varepsilon_j=\pm1
}\left\|\sum\limits_{j=0}^n \varepsilon_j x_j\right\|,$$
where $\{r_j(t)\}_{j=0}^n$, $n\in\mathbb{N}$ is a set of Rademacher functions, i.e. $r_j(t)=sign(\sin(2^{j+1} \pi t))$, $j=0,1,2\dots$. Introduce the following definitions.
\begin{definition}\label{OR-type}
We will say that a Banach space $E$ has Orlicz-Rademacher type $\Phi$ with constant $T_{\Phi}(E)$ provided there exists an Orlicz function $\Phi$ and a constant $T_{\Phi}(E)$ so that, for every finite set of vectors $\{x_j\}_{j=0}^{n}\subset E$, we have
\begin{equation}\label{O-type}
    \left(\mathbb{E}\left\|\sum\limits_{j=0}^{n} \varepsilon_j x_j \right\|^2\right)^{\frac{1}{2}}
\leq T_{\Phi}(E) \inf\left\{\rho >0:\quad  \sum\limits_{j=0}^{n} \Phi\left(\frac{\|x_j\|}{\rho}\right)\leq 1\right\}
\end{equation}
\end{definition}
\begin{definition}\label{OR-cotype}
We will say that a Banach space $E$ has Orlicz-Rademacher cotype $\Psi$ with constant $C_{\Psi}(E)$ provided there exists an Orlicz function $\Psi$ and a constant $C_{\Psi}(E)$ so that, for every finite set of vectors $\{x_j\}_{j=0}^{n}\subset E$, we have
\begin{equation}\label{O-cotype}
\left(\mathbb{E}\left\|\sum\limits_{j=0}^{n} \varepsilon_j x_j \right\|^2\right)^{\frac{1}{2}}
\geq C_{\Psi}(E) \inf\left\{\rho >0:\quad  \sum\limits_{j=0}^{n} \Psi\left(\frac{\|x_j\|}{\rho}\right)\leq 1\right\}
\end{equation}
\end{definition}
These notions are natural generalizations of the well known concepts of type and cotype of a Banach space $E$.
Clearly, if we take $\Phi(t)=t^p$ with $p\geq 1$, then $\inf\left\{\rho >0:\quad  \sum\limits_{j=0}^{n} \Phi\left(\frac{\|x_j\|}{\rho}\right)\leq 1\right\}=\left(\sum\limits_{j=0}^{n} \|x_n\|^p\right)^{\frac{1}{p}}$ and Definition \ref{OR-type} coincides with definition of Banach space of (Rademacher) type $p$ \cite{Johnson2}. Analogously, if we take $\Psi(t)=t^q$ with $q\geq 1$, Definition \ref{OR-cotype} will coincide with definition of Banach space of (Rademacher) cotype $q$ \cite{Johnson2}. In fact, all the possible values of the type of a Banach space are located on $[1,2]$ and all the possible values of the cotype of a Banach space belong to $[2,\infty]$. Every Banach space $E$ has a trivial type $p=1$ and trivial cotype $q=\infty$ $\left( \text{with}\:
\max\limits_{0\leq j\leq n} \|x_j\|\: \text{replacing the expression}\:\left(\sum\limits_{j=0}^{n} \|x_j\|^q\right)^{\frac{1}{q}} \:\text{when}\: q=\infty\right)$\newline by convexity of the norm \cite{Johnson2,Lindenstrauss}. By $L_p(\mu)$, $p\in [1,\infty)$, we mean a Banach space of $\mu$-measurable functions $f$ for which $\|f\|=\left(\int |f|^p d\mu\right)^{\frac{1}{p}}<\infty$. The spaces $L_p(\mu)$, $p\in [1,\infty)$, have (the best possible) type $\min\{2,p\}$ and (the best possible) cotype $\max\{2,p\}$ \cite{Johnson2,Lindenstrauss}. As in the definitions of type and cotype, the $L_2$ average $\left(\int\limits_{0}^{1} \left\|\sum\limits_{j=0}^{n} r_j (t) x_j
\right\|^2 dt\right)^{\frac{1}{2}}$ can be replaced in Definitions \ref{OR-type}, \ref{OR-cotype} by any other $L_p$ average, $p\in[1,\infty)\setminus {2}$, without affecting the definition. But of course then the constants in (\ref{O-type}), (\ref{O-cotype}) will change. This follows from the Kahane-Khintchine inequality \cite{Johnson2, Lindenstrauss}. Along with Definitions  \ref{OR-type}, \ref{OR-cotype} we intoduce the following.
\begin{definition}\label{OR-infratype}
We will say that a Banach space $E$ has Orlicz-Rademacher infratype $\Phi$ with constant $I_{\Phi}(E)$ provided there exists an Orlicz function $\Phi$ and a constant $I_{\Phi}(E)$ so that, for every finite set of vectors $\{x_j\}_{j=0}^{n}\subset E$, we have
\begin{equation}\label{O-infratype}
\min_{\varepsilon_j=\pm 1} \left\| \sum\limits_{j=0}^n \varepsilon_j x_j  \right\|
\leq I_{\Phi}(E) \inf\left\{\rho >0:\quad  \sum\limits_{j=0}^{n} \Phi\left(\frac{\|x_j\|}{\rho}\right)\leq 1\right\}
\end{equation}
\end{definition}
\begin{definition}\label{OR-M-cotype}
We will say that a Banach space $E$ has Orlicz-Rademacher $M$-cotype $\Psi$ with constant $M_{\Psi}(E)$, if there exists an Orlicz function $\Psi$ and a constant $M_{\Psi}(E)$ so that, for every finite set of vectors $\{x_j\}_{j=0}^{n}\subset E$, we have
\begin{equation}\label{O-M-cotype}
\max_{\varepsilon_j=\pm 1} \left\| \sum\limits_{j=0}^n \varepsilon_j x_j  \right\|
\geq M_{\Psi}(E) \inf\left\{\rho >0:\quad  \sum\limits_{j=0}^{n} \Psi\left(\frac{\|x_j\|}{\rho}\right)\leq 1\right\}
\end{equation}
\end{definition}
The latter two concepts may be considered as natural generalizations of the concepts of infratype and $M$-cotype of a Banach space $E$ \cite{Kadets}. This can be done in a way similar to above mentioned. Infratype and $M$-cotype have properties that are similar to the properties of type and cotype. E.g., if an infinite-dimensional Banach space $E$ has infratype $p$ and $M$-cotype $q$, then $p\leq 2$ and $q\geq 2$ \cite{Kadets}. The spaces $L_p(\mu)$, $p\in [1,\infty)$, have (the best possible) infratype $\min\{2,p\}$ and (the best possible) $M$-cotype $\max\{2,p\}$ \cite{Kadets}.
One can find more about type, infratype, cotype and $M$-cotype in \cite{Johnson2,Johnson3,Lindenstrauss,Kadets}.
\subsection{Unconditional Schauder decompositions and Banach space geometry}
The following Lemma shows that the properties of unconditional Schauder decompositions in $E$ depend directly on the character of intrinsic geometric structure of $E$.
\begin{lemma}\label{Lemma}
Let $E$ be a Banach space which possesses unconditional\newline Schauder decomposition $\{\mathfrak{M}_n\}_{n=0}^{\infty}$ with constant $M$ and the a.s.c.p. $\{P_n\}_{n=0}^{\infty}$. Also assume that, on the one hand, $E$ has or Orlicz-Rademacher type $\Phi$, or Orlicz-Rademacher infratype $\Phi$, and, on the other hand, $E$ has or Orlicz-Rademacher cotype $\Psi$, or Orlicz-Rademacher $M$-cotype $\Psi$.

Then there exist constants $T=T(\Phi,M)>0$ and $C=C(\Psi,M)>0$ such that for each $x\in E$ we have
$$C \inf\left\{\rho >0:\quad  \sum\limits_{n=0}^{\infty} \Psi\left(\frac{\|P_n x\|}{\rho}\right)\leq 1\right\} \leq \|x\|$$
\begin{equation}\label{Type-Cotype}
\leq T  \inf\left\{\rho >0:\quad  \sum\limits_{n=0}^{\infty} \Phi\left(\frac{\|P_n x\|}{\rho}\right)\leq 1\right\}.
\end{equation}
\end{lemma}
In particular, for unconditional Schauder decompositions in Banach spaces with classical geometric characteristics, such as classical type, infratype, cotype and $M$-cotype we have the following simple consequence of Lemma \ref{Lemma}.
\begin{corollary}\label{Corollary}
Let $E$ be a Banach space which possesses unconditional\newline Schauder decomposition $\{\mathfrak{M}_n\}_{n=0}^{\infty}$ with constant $M$ and the a.s.c.p. $\{P_n\}_{n=0}^{\infty}$. Also assume that, on the one hand, $E$ has or type $p$, or infratype $p$, and, on the other hand, $E$ has or cotype $q$, or $M$-cotype $q$.

Then there exist constants $T=T(p,M)>0$ and $C=C(q,M)>0$ such that for each $x\in E$ we have
\begin{equation}\label{Typep-Cotypeq}
C \left(\sum\limits_{n=0}^{\infty} \|P_n x\|^q\right)^{\frac{1}{q}}  \leq \|x\|\leq T  \left(\sum\limits_{n=0}^{\infty} \|P_n x\|^p\right)^{\frac{1}{p}}.
\end{equation}
\end{corollary}
\begin{remark}\label{p=q=2}
Let $\{\mathfrak{M}_n\}_{n=0}^{\infty}$ be an unconditional Schauder decomposition of Hilbert space $H$ with constant $M$ and the a.s.c.p. $\{P_n\}_{n=0}^{\infty}$. Then, for each $x\in H$ one has the following inequality, which coincides with (\ref{2}).
\begin{equation}\label{p=q=2 inequality}
\frac{1}{2M} \left(\sum\limits_{n=0}^{\infty} \|P_n x\|^2\right)^{\frac{1}{2}}  \leq \|x\|\leq 2M  \left(\sum\limits_{n=0}^{\infty} \|P_n x\|^2\right)^{\frac{1}{2}}.
\end{equation}
\end{remark}
\begin{remark}\label{Ball}
Under the assumptions of Lemma \ref{Lemma}, for every $x\in E$ one has the following inequality
$$C \|\:\{\|P_n x\|\}_{n=0}^{\infty}\:\|_{\ell_{\Psi}} \leq \|x\|\leq T  \|\:\{\|P_n x\|\}_{n=0}^{\infty}\:\|_{\ell_{\Phi}},$$
which reads that each $x\in E$ is located, on the one hand, inside the ball $B_{\Phi}$ with center $0$ and radius $T  \|\:\{\|P_n x\|\}_{n=0}^{\infty}\:\|_{\ell_{\Phi}}$ and, on the other hand, outside of the ball $B_{\Psi}$ with center $0$ and radius  $C \|\:\{\|P_n x\|\}_{n=0}^{\infty}\:\|_{\ell_{\Psi}}$.
\end{remark}
To demonstrate how does this theory work for the case of concrete Banach spaces  with unconditional Schauder decompositions we consider, for example, unconditional Schauder decompositions in $L_p(\mu)$ spaces.
For the case of $L_p(\mu)$ spaces it is known that the best constants in inequalities (\ref{O-type}), (\ref{O-cotype}), which correspond to the best possible type $\min\{2,p\}$ and the best possible cotype $\max\{2,p\}$, satisfy the following relations
$$T_p(L_p(\mu))=1,\:\: C_2(L_p(\mu))\geq A_p\:\:\text{for}\:\: p\in[1,2],$$
$$T_2(L_p(\mu))\leq B_p,\:\: C_p(L_p(\mu))=1\:\:\text{for}\:\: p\in[2,\infty),$$
where $A_p$ and $B_p$ are the best constants in the Khintchine inequality \cite{Johnson2}. As it was shown by Haagerup in \cite{Haagerup},
$$A_p=\left \{
\begin{array}{l}
\displaystyle 2^{\frac{1}{2}-\frac{1}{p}},\quad\quad\quad\quad\quad\quad p\in (0,p_0],\\[10pt]\displaystyle
\displaystyle 2^{\frac{1}{2}}\left(\frac{\Gamma \left( \frac{p+1}{2} \right)}{\sqrt{\pi}} \right)^{\frac{1}{p}},\quad p\in [p_0,2],\\[10pt]\displaystyle
1,\quad\quad\quad\quad\quad\quad\quad\quad p\in[2,\infty),\\
\end{array}
\right .$$
$$B_p=\left \{
\begin{array}{l}
\displaystyle 1,\quad\quad\quad\quad\quad\quad\quad\quad p\in(0,2],\\[10pt]\displaystyle
\displaystyle 2^{\frac{1}{2}}\left(\frac{\Gamma \left( \frac{p+1}{2} \right)}{\sqrt{\pi}} \right)^{\frac{1}{p}},\quad p\in [2,\infty),\\
\end{array}
\right .$$
where $p_0$ is the solution of the equation $\Gamma  \left(\frac{p+1}{2}\right)=\frac{\sqrt{\pi}}{2}$ in the interval $[1,2]$, $p_0\approx 1,84742$, and $\Gamma$ is the gamma function. Consequently, for $p\in [1,p_0]$ we have $C_2(L_p(\mu))\geq 2^{\frac{1}{2}-\frac{1}{p}}$, for $p\in [p_0,2]$ we have $C_2(L_p(\mu))\geq 2^{\frac{1}{2}}\left(\frac{\Gamma \left( \frac{p+1}{2} \right)}{\sqrt{\pi}} \right)^{\frac{1}{p}}$, and for $p\in [2,\infty)$ we have that $T_2(L_p(\mu))\leq 2^{\frac{1}{2}}\left(\frac{\Gamma \left( \frac{p+1}{2} \right)}{\sqrt{\pi}} \right)^{\frac{1}{p}}$. Taking these considerations into account, the application of Corollary \ref{Corollary} leads to the following proposition on unconditional Schauder decompositions in $L_p(\mu)$ spaces.
\begin{proposition}\label{Proposition}
Let $\{\mathfrak{M}_n\}_{n=0}^{\infty}$ be an unconditional Schauder decomposition of $L_p(\mu)$ with constant $M$ and the a.s.c.p. $\{P_n\}_{n=0}^{\infty}$.

(i) If $p\in [1,p_0]$, then for each $x\in L_p(\mu)$ we have
\begin{align}
\label{1<p<p_0 inequality}
\frac{2^{-\frac{1}{2}-\frac{1}{p}}}{M} \left(\sum\limits_{n=0}^{\infty} \|P_n x\|^2\right)^{\frac{1}{2}} \leq \|x\| &\leq 2M   \left(\sum\limits_{n=0}^{\infty} \|P_n x\|^p\right)^{\frac{1}{p}}.\\
\intertext{(ii) If $p\in [p_0,2]$, then for each $x\in L_p(\mu)$ we have}
\label{p_0<p<2 inequality}
\frac{1}{\sqrt{2} M}  \left(\frac{\Gamma \left( \frac{p+1}{2} \right)}{\sqrt{\pi}} \right)^{\frac{1}{p}}  \left(\sum\limits_{n=0}^{\infty} \|P_n x\|^2\right)^{\frac{1}{2}} \leq \|x\| &\leq 2M   \left(\sum\limits_{n=0}^{\infty} \|P_n x\|^p\right)^{\frac{1}{p}}.
\end{align}
\end{proposition}
(iii) If $p\in [2,\infty)$, then for each $x\in L_p(\mu)$ we have
$$\frac{1}{2 M}   \left(\sum\limits_{n=0}^{\infty} \|P_n x\|^p\right)^{\frac{1}{p}} \leq \|x\|$$
\begin{equation}\label{2<p<infty inequality}
    \leq  \sqrt{8} \left(\frac{\Gamma \left( \frac{p+1}{2} \right)}{\sqrt{\pi}} \right)^{\frac{1}{p}} M \left(\sum\limits_{n=0}^{\infty} \|P_n x\|^2\right)^{\frac{1}{2}}.
\end{equation}

This assertion will be very useful in the derivation of the stability theorem for unconditional Schauder decompositions in $\ell_p$ spaces, see Section 5, Theorem \ref{ell_p Theorem}.
\subsection{Isomorphic unconditional Schauder decompositions}
In what follows we denote by $\{e_n\}_{n=0}^{\infty}$ the canonical basis of $\ell_p, \:p\in [1,\infty)$, i.e. $e_n=(\delta_i^n)$, $n\in\mathbb{Z}_{+}$. Let $h_{\Phi}$ be a subset of Orlicz sequence space $\ell_{\Phi}$ consisting of sequences $x=(a_0,a_1,\dots)\in \ell_{\Phi}$ such that  $\sum\limits_{n=0}^{\infty} \Phi\left(\frac{|a_n|}{\rho}\right)<\infty$ for every $\rho>0$. Clearly, $h_{\Phi}$ is a closed subspace of  $\ell_{\Phi}$ and $\{e_n\}_{n=0}^{\infty}$ forms an unconditional basis of $h_{\Phi}$. The properties of spaces $\ell_{\Phi}$ are inseparably connected with the properties of Orlicz functions which generate them. E.g., if $\Phi(t)=0$ for some $t>0$, then $\ell_{\Phi}$ is isomorphic to $\ell_{\infty}$, and $h_{\Phi}$ is isomorphic to $c_0$ \cite{Lindenstrauss}. Hence, in general, the spaces $\ell_{\Phi}$ and $h_{\Phi}$ are distinct. In order to impose conditions for $\ell_{\Phi}$ to coincide with $h_{\Phi}$ we need the following.
\begin{definition}[\cite{Lindenstrauss}]\label{2-condition}
An Orlicz function $\Phi$ is said to satisfy the $\Delta_2$-condition at zero provided $\limsup\limits_{t\rightarrow 0} \frac{\Phi(2t)}{\Phi(t)}<\infty$.
\end{definition}
For an Orlicz function $\Phi$ we have the following result.
\begin{proposition}[\cite{Lindenstrauss}]\label{M-proposition}
Let $\Phi$ be an Orlicz function. Then the following statements are equivalent.

(1) $\Phi$ satisfies the $\Delta_2$-condition at zero.

(2) $\ell_{\Phi}=h_{\Phi}$.

(3) $\ell_{\Phi}$ is separable.

(4) $\ell_{\Phi}$ contains no subspace isomorphic to $\ell_{\infty}$.
\end{proposition}
Let an Orlicz function $\Phi$ satisfies the $\Delta_2$-condition at zero. Then, by virtue of Proposition \ref{M-proposition}, $\ell_{\Phi}=h_{\Phi}$ and, hence, $\{e_n\}_{n=0}^{\infty}$ forms an unconditional basis of $\ell_{\Phi}$. In order to give some stability properties for arbitrary unconditional Schauder decompositions in $E$ we need the following.
\begin{definition}[\cite{Singer2}]\label{ellphi-Besselian}
A Schauder decomposition $\{\mathfrak{M}_n\}_{n=0}^{\infty}$ of a Banach space $E$ is called $\ell_{\Phi}$-Besselian ($\infty$-Besselian), if the convergence of $\sum\limits_{n=0}^{\infty} x_n$ in $E$, where $x_n\in \mathfrak{M}_n$ $n=0,1,\dots$, implies the convergence of $\sum\limits_{n=0}^{\infty} \|x_n\| e_n$ in $\ell_{\Phi}$ ($c_0$).
\end{definition}
\begin{definition}[\cite{Singer2}]\label{ellphi-Hilbertian}
A Schauder decomposition $\{\mathfrak{M}_n\}_{n=0}^{\infty}$ of a Banach space $E$ is called $\ell_{\Phi}$-Hilbertian ($\infty$-Hilbertian), if the convergence of \newline$\sum\limits_{n=0}^{\infty} \|x_n\| e_n$ in $\ell_{\Phi}$ ($c_0$)  implies the convergence of $\sum\limits_{n=0}^{\infty} x_n$ in $E$, where $x_n\in \mathfrak{M}_n$ $n=0,1,\dots$.
\end{definition}
In the particular case when $\ell_{\Phi}=\ell_p$, $p\in[1,\infty)$, we use the terms $p$-Besselian, $p$-Hilbertian Schauder decompositions.
A Schauder decomposition which is both $p$-Besselian and $p$-Hilbertian is called an $\ell_p$-decomposition if $p<\infty$ and a $c_0$-decomposition if $p=\infty$ \cite{Singer2}. Obviously, each Schauder decomposition $\{\mathfrak{M}_n\}_{n=0}^{\infty}$ of $E$ is both $\infty$-Besselian and $1$-Hilbertian. We have the following characterization of $\ell_{\Phi}$-Besselian and $\ell_{\Phi}$-Hilbertian Schauder decompositions.
\begin{theorem}[\cite{Singer2}]\label{ellPhi-Besselian}
Let $\{\mathfrak{M}_n\}_{n=0}^{\infty}$ be a Schauder decomposition of a Banach space $E$. Then

(i) $\{\mathfrak{M}_n\}_{n=0}^{\infty}$ is $\ell_{\Phi}$-Besselian ($\infty$-Besselian) if and only if there exists a constant $c>0$ such that
$$c\left\|\sum\limits_{j=0}^{n} \|x_j\| e_j\right\|_{\ell_{\Phi}(c_0)}\leq \left\|\sum\limits_{j=0}^{n} x_j\right\|$$
is satisfied for all finite sequences $x_j\in\mathfrak{M}_j$ $j=0,\dots,n.$

(ii) $\{\mathfrak{M}_n\}_{n=0}^{\infty}$ is $\ell_{\Phi}$-Hilbertian ($\infty$-Hilbertian) if and only if there exists a constant $C>0$ such that
$$C\left\|\sum\limits_{j=0}^{n} \|x_j\| e_j\right\|_{\ell_{\Phi}(c_0)}\geq \left\|\sum\limits_{j=0}^{n} x_j\right\|$$
is satisfied for all finite sequences $x_j\in\mathfrak{M}_j$ $j=0,\dots,n.$
\end{theorem}
\begin{remark}\label{ell_2-decomposition}
Combining Theorem \ref{ellPhi-Besselian} with Theorem \ref{Orthogonal}, we obtain that a Schauder decomposition $\{\mathfrak{M}_n\}_{n=0}^{\infty}$ in Hilbert space $H$ is unconditional if and only if $\{\mathfrak{M}_n\}_{n=0}^{\infty}$ is an $\ell_2$-decomposition.
\end{remark}
A similar situation with $\ell_p(c_0)$-decompositions and unconditional decompositions for the case of Banach spaces is the exception rather than the rule. Combining Theorem \ref{ellPhi-Besselian} with Theorem \ref{USD}, we clearly see that every $\ell_p$-decomposition or every $c_0$-decomposition is unconditional, but the converse, in general, is false. Among the exceptions we specify the spaces $\ell_1$ and $c_0$. Namely, applying of Theorem 15.19 of \cite{Singer2} yields the following.
\begin{proposition}\label{ell_1,c_0}
(i) A Schauder decomposition $\{\mathfrak{M}_n\}_{n=0}^{\infty}$ in $\ell_1$ is unconditional if and only if $\{\mathfrak{M}_n\}_{n=0}^{\infty}$ is an $\ell_1$-decomposition.
(ii) A Schauder decomposition $\{\mathfrak{M}_n\}_{n=0}^{\infty}$ in $c_0$ is unconditional if and only if $\{\mathfrak{M}_n\}_{n=0}^{\infty}$ is a $c_0$-decomposition.
\end{proposition}
This situation is the reminiscent of the situation with the uniqueness of unconditional basis in Banach space. In this context, the question of the uniqueness of unconditional Schauder decomposition in space $\ell_1$ or in space $c_0$ becomes natural. On the other hand, for any unconditional Schauder decomposition in space with certain Orlicz-Rademacher structural properties we have, by virtue of Lemma \ref{Lemma} and Theorem \ref{ellPhi-Besselian}, the following result.
\begin{theorem}\label{ThBesselian}
Let $E$ be a Banach space with unconditional Schauder decomposition $\{\mathfrak{M}_n\}_{n=0}^{\infty}$. Further, assume that, on the one hand, $E$ has or Orlicz-Rademacher type $\Phi$, or Orlicz-Rademacher infratype $\Phi$, and, on the other hand, $E$ has or Orlicz-Rademacher cotype $\Psi$, or Orlicz-Rademacher $M$-cotype $\Psi$.

Then $\{\mathfrak{M}_n\}_{n=0}^{\infty}$ is both $\ell_{\Phi}$-Hilbertian and $\ell_{\Psi}$-Besselian.
\end{theorem}
For any unconditional Schauder decomposition in space with classical type, infratype, cotype and $M$-cotype we have the following simple consequence of Theorem \ref{ThBesselian}.
\begin{corollary}\label{Th p-Besselian}
Let $E$ be a Banach space with unconditional Schauder decomposition $\{\mathfrak{M}_n\}_{n=0}^{\infty}$. Further, assume that, on the one hand, $E$ has or type p, or infratype p, and, on the other hand, $E$ has or cotype q, or $M$-cotype q.

Then $\{\mathfrak{M}_n\}_{n=0}^{\infty}$ is both p-Hilbertian and q-Besselian.
\end{corollary}
In order to formulate stability theorem for unconditional Schauder decompositions in $E$ we need the following.
\begin{definition}[\cite{Singer2}]\label{biorthogonal}
A pair of sequences $(\{\mathfrak{M}_n\}_{n=0}^{\infty},\{P_n\}_{n=0}^{\infty})$, where\newline $\{\mathfrak{M}_n\}_{n=0}^{\infty}$ is a sequence of closed nonzero linear subspaces of a Banach space $E$ and $\{P_n\}_{n=0}^{\infty}$ is a sequence of bounded linear projections satisfying $P_n E=\mathfrak{M}_n$ for all $n$, will be called a generalized biorthogonal system provided it satisfies $P_i P_j=\delta_{i}^{j} P_i= \delta_{i}^{j} P_j$ for $i,j=0,1,\dots$ The generalized biorthogonal system $(\{\mathfrak{M}_n\}_{n=0}^{\infty},\{P_n\}_{n=0}^{\infty})$ is said to be $E$-complete, if $\overline{Lin} \{\mathfrak{M}_n\}_{n=0}^{\infty} = E$.
\end{definition}
Obviously, if $\{\mathfrak{M}_n\}_{n=0}^{\infty}$ is a Schauder decomposition of $E$ with the a.s.c.p. $\{P_n\}_{n=0}^{\infty}$, then $(\{\mathfrak{M}_n\}_{n=0}^{\infty},\{P_n\}_{n=0}^{\infty})$ is an $E$-complete generalized biorthogonal system, but the converse is not true.
\begin{definition}[\cite{Singer2}]\label{omega}
A sequence of nonzero subspaces of $E$ is said to be $\omega$-linearly independent, if the relations $\sum\limits_{n=0}^{\infty} x_n =0$, $x_n\in \mathfrak{M}_n$ $(n=0,1,\dots)$, imply $x_n=0$ $(n=0,1,\dots)$.
\end{definition}
Combining Corollary \ref{Th p-Besselian} with Theorem 15.17 of \cite{Singer2}, we obtain the following stability theorem which is valid for every unconditional Schauder decomposition in a Banach space $E$.
\begin{theorem}\label{IsomorphicUSD}
Let $E$ has or cotype $q$, or $M$-cotype $q$, where $q\in[2,\infty]$, and assume that $E$ has an unconditional Schauder decomposition $\{\mathfrak{M}_n\}_{n=0}^{\infty}$ with the a.s.c.p. $\{P_n\}_{n=0}^{\infty}$. Then

(i) There exists a constant $\lambda\in(0,1)$ such that every sequence of subspaces $\{\mathfrak{N}_n\}_{n=0}^{\infty}$ of $E$, satisfying
\begin{equation}\label{openings}
    \left(\sum\limits_{n=0}^{\infty} \theta\left(\mathfrak{M}_n,\mathfrak{N}_n\right)^p  \right)^{\frac{1}{p}} \leq \lambda,
\end{equation}
where $p^{-1}+q^{-1}=1$ (if $q=\infty$, then $p=1$) and\newline $\theta\left(\mathfrak{M},\mathfrak{N}\right)=\max \left\{\sup\limits_{x\in \mathfrak{M}, \|x\|=1} dist (x,\mathfrak{N}),\sup\limits_{y\in \mathfrak{N}, \|y\|=1} dist (y,\mathfrak{M}) \right\}$ is the opening of the subspaces $\mathfrak{M}$, $\mathfrak{N}$, is a Schauder decomposition of $E$, isomorphic to $\{\mathfrak{M}_n\}_{n=0}^{\infty}$.
Moreover, a constant $\lambda$ may be chosen in the following way
$$\lambda=\frac{1}{4 \sup\limits_{0\leq n<\infty} \left\|\sum\limits_{j=0}^{n} P_j \right\| \left(1+ \sup\limits_{0\leq n<\infty} \|P_n\| \right)^2}.$$

(ii) Every sequence of subspaces $\{\mathfrak{N}_n\}_{n=0}^{\infty}$ of $E$, satisfying
\begin{equation}\label{serie of openings}
    \sum\limits_{n=0}^{\infty} \theta\left(\mathfrak{M}_n,\mathfrak{N}_n\right)^p  < \infty,
\end{equation}
where $p^{-1}+q^{-1}=1$ ($p=1$, if $q=\infty$), and admitting a sequence $\{J_n\}_{n=0}^{\infty}$ such that $(\{\mathfrak{N}_n\}_{n=0}^{\infty},\{J_n\}_{n=0}^{\infty})$ is an $E$-complete generalized biorthogonal system, is a $q$-Besselian Schauder decomposition of $E$. If, additionally, $\{\mathfrak{M}_n\}_{n=0}^{\infty}$ is an FDD, then the same conclusion holds for every $\omega$-linearly independent sequence of subspaces $\{\mathfrak{N}_n\}_{n=0}^{\infty}$ satisfying (\ref{serie of openings}).
\end{theorem}
\begin{remark}\label{IsomorphicUSDremark}
Note that, by virtue of Theorem \ref{USD}, every sequence of subspaces $\{\mathfrak{N}_n\}_{n=0}^{\infty}$, isomorphic to unconditional Schauder decomposition\newline $\{\mathfrak{M}_n\}_{n=0}^{\infty}$ with constant $M$, is itself an unconditional Schauder decomposition with constant $M \|S\|\|S^{-1}\|$, where $\mathfrak{N}_n=S\mathfrak{M}_n$, $n=0,1,\dots$
\end{remark}

\section{The main results}
\subsection{Isomorphic Schauder decompositions in Orlicz spaces $\ell_{\Phi}$}
Passing to the limit in Theorem \ref{ellPhi-Besselian} yields that a Schauder decomposition $\{\mathfrak{M}_n\}_{n=0}^{\infty}$ of $E$ with the a.s.c.p. $\{P_n\}_{n=0}^{\infty}$ is $\ell_{\Psi}$-Hilbertian ($\infty$-Hilbertian) if and only if there exists a constant $C>0$ such that for every $x\in E$
\begin{equation}\label{withC}
    \left\|x\right\| \leq C\left\|\sum\limits_{n=0}^{\infty} \|P_n x\| e_n\right\|_{\ell_{\Psi}(c_0)}.
\end{equation}
Thus, we arrive at the definition.
\begin{definition}\label{with constant C}
We will say that a Schauder decomposition $\{\mathfrak{M}_n\}_{n=0}^{\infty}$ of $E$ is $\ell_{\Psi}$-Hilbertian ($\infty$-Hilbertian) with constant $C$ if and only if (\ref{withC}) holds for each $x\in E$.
\end{definition}
We have the following stability result for $\ell_{\Psi}$-Hilbertian ($\infty$-Hilbertian) \newline Schauder decompositions $\{\mathfrak{M}_n\}_{n=0}^{\infty}$ in Orlicz spaces $\ell_{\Phi}$.
\begin{theorem}\label{MainTheorem1}
Let $\{\mathfrak{M}_n\}_{n=0}^{\infty}$ be an $\ell_{\Psi}$-Hilbertian ($\infty$-Hilbertian) Schauder decomposition of the Orlicz space $\ell_{\Phi}$ with constant $C$ and the a.s.c.p. $\{P_n\}_{n=0}^{\infty}$, such that for every $x=(a_0,a_1,a_2,\dots)\in \ell_{\Phi}$ one has $$P_0 x=(a_0,a_1,a_2,\dots,a_{m-1},0,0,\dots).$$ Assume that $\{J_n\}_{n=0}^{\infty}$ is a sequence of nonzero projections in $\ell_{\Phi}$ satisfying $J_n J_m=\delta_{n}^{m} J_n$ for $n,m=0,1,\dots$ Furthermore, suppose that the condition (\ref{1}) holds and for all $x\in \ell_{\Phi}$ we have
\begin{equation}\label{*}
    \left\|\sum\limits_{n=1}^{\infty} \|P_n (J_n-P_n) x\| e_n \right\|_{\ell_{\Psi}(c_0)} \leq \varsigma \|x\|,
\end{equation}
where $\varsigma$ is a constant such that $0\leq \varsigma<C^{-1}$.

Then $\{\mathfrak{M}_n\}_{n=0}^{\infty}$ and $\{J_n(\ell_{\Phi})\}_{n=0}^{\infty}$ are isomorphic Schauder decompositions of $\ell_{\Phi}$.
\end{theorem}
Since every Schauder decomposition is $1$-Hilbertian with constant $1$, we obtain the following consequence of Theorem \ref{MainTheorem1} concerning stability of Schauder decompositions of certain structure in Orlicz spaces $\ell_{\Phi}$.
\begin{corollary}\label{1-Hilbertian}
Let $\{\mathfrak{M}_n\}_{n=0}^{\infty}$ be a Schauder decomposition of $\ell_{\Phi}$ with the a.s.c.p. $\{P_n\}_{n=0}^{\infty}$, such that for every \newline $x=(a_0,a_1,a_2,\dots)\in \ell_{\Phi}$ one has $$P_0 x=(a_0,a_1,a_2,\dots,a_{m-1},0,0,\dots).$$ Assume that $\{J_n\}_{n=0}^{\infty}$ is a sequence of nonzero projections in $\ell_{\Phi}$ satisfying $J_n J_m=\delta_{n}^{m} J_n$ for $n,m=0,1,\dots$ Furthermore, suppose that the condition (\ref{1}) holds and for all $x\in \ell_{\Phi}$ we have
$$\sum\limits_{n=1}^{\infty} \|P_n (J_n-P_n) x\|\leq c \|x\|,$$
where $c$ is a constant such that $0\leq c<1$. Then $\{\mathfrak{M}_n\}_{n=0}^{\infty}$ and $\{J_n(\ell_{\Phi})\}_{n=0}^{\infty}$ are isomorphic Schauder decompositions of $\ell_{\Phi}$.
\end{corollary}
In case when all the subspaces of $\{\mathfrak{M}_n\}_{n=0}^{\infty}$ are one dimensional, by Corollary \ref{1-Hilbertian}, we have the following result on stability of bases of certain structure in $\ell_{\Phi}$ spaces.
\begin{corollary}\label{1-Hilbertian bases}
Let $\{\phi_n\}_{n=0}^{\infty}$ be a basis of $\ell_{\Phi}$ with corresponding sequence of coordinate functionals $\{\phi_n^{\ast}\}_{n=0}^{\infty}$, such that for all $x=(a_0,a_1,a_2,\dots)\in \ell_{\Phi}$ one has
$$\langle\phi_0^{\ast},x\rangle \phi_0=(a_0,0,0,\dots).$$
Let $(\{\psi_n\}_{n=0}^{\infty},\{\psi_n^{\ast}\}_{n=0}^{\infty})$ be a biorthogonal system in $\ell_{\Phi}$, i.e. $\langle\psi_i^{\ast},\psi_j\rangle =\delta_{i}^{j}$ for $i,j=0,1,2,\dots$. If for all $x\in\ell_{\Phi}$ we have
$$\sum\limits_{n=1}^{\infty} |\langle\psi_n^{\ast}, x\rangle \langle\phi_n^{\ast}, \psi_n\rangle - \langle\phi_n^{\ast}, x\rangle  | \|\phi_n\|\leq c \|x\|,$$
where $c$ is a constant such that $0\leq c<1$, then $\{\phi_n\}_{n=0}^{\infty}$ and $\{\psi_n\}_{n=0}^{\infty}$ are isomorphic bases of the space $\ell_{\Phi}$.
\end{corollary}
To demonstrate what does Theorem \ref{MainTheorem1} say for the case of concrete Orlicz spaces possessing Schauder decompositions we consider, for instance, unconditional Schauder decompositions in $\ell_p$ spaces.
Combining Theorem \ref{MainTheorem1} with Proposition \ref{Proposition}, in view of Remark \ref{IsomorphicUSDremark}, we obtain the following stability result for unconditional Schauder decompositions in $\ell_p$ spaces.
\begin{theorem}\label{ell_p Theorem}
Let $\{\mathfrak{M}_n\}_{n=0}^{\infty}$ be an unconditional Schauder decomposition of the space $\ell_p$, $p\in[1,\infty)$, with constant $M$ and the a.s.c.p. $\{P_n\}_{n=0}^{\infty}$, such that for every $x=(a_0,a_1,a_2,\dots)\in \ell_p$ one has $$P_0 x=(a_0,a_1,a_2,\dots,a_{m-1},0,0,\dots).$$  Assume that $\{J_n\}_{n=0}^{\infty}$ is a sequence of nonzero projections in $\ell_p$ satisfying (\ref{1}), such that $J_n J_m=\delta_{n}^{m} J_n$ for $n,m=0,1,\dots$ Then the following statements hold

(i) If $p\in [1,2]$ and for every $x\in \ell_p$ we have
$$\left(\sum\limits_{n=1}^{\infty} \|P_n (J_n-P_n) x\|^p\right)^{\frac{1}{p}}\leq \varsigma_1 \|x\|,$$
where $\varsigma_1\in \left[0,(2M)^{-1}\right)$,
then $\{J_n(\ell_p)\}_{n=0}^{\infty}$ is also an unconditional Schauder decomposition of $\ell_p$, $p\in [1,2]$, isomorphic to $\{\mathfrak{M}_n\}_{n=0}^{\infty}$.

(ii) If $p\in [2,\infty)$ and for every $x\in \ell_p$ we have
$$\left(\sum\limits_{n=1}^{\infty} \|P_n (J_n-P_n) x\|^2\right)^{\frac{1}{2}}\leq \varsigma_2 \|x\|,$$
where $\varsigma_2\in \left[0,\frac{1}{\sqrt{8} M} \left(\frac{\Gamma \left( \frac{p+1}{2} \right)}{\sqrt{\pi}} \right)^{-\frac{1}{p}} \right)$,
then $\{J_n(\ell_p)\}_{n=0}^{\infty}$ is also an unconditional Schauder decomposition of $\ell_p$, $p\in [2,\infty)$, isomorphic to $\{\mathfrak{M}_n\}_{n=0}^{\infty}$.
\end{theorem}
We say that a sequence of elements $\{\phi_n\}_{n=0}^{\infty}$ is unconditional basis of $E$ with constant $M$, provided that the sequence of corresponding one dimensional subspaces $\{Lin \{ \phi_n\}\}_{n=0}^{\infty}$ forms an unconditional Schauder decomposition of $E$ with constant $M$.
In case when all the subspaces of $\{\mathfrak{M}_n\}_{n=0}^{\infty}$ are one dimensional, by virtue of Theorem \ref{ell_p Theorem}, we have the following result on stability of unconditional bases in $\ell_p$ spaces.
\begin{corollary}\label{ell p bases}
Let $\{\phi_n\}_{n=0}^{\infty}$ be an unconditional basis of $\ell_p$, $p\in[1,\infty)$, with constant $M$ and corresponding sequence of coordinate functionals $\{\phi_n^{\ast}\}_{n=0}^{\infty}$, such that for all $x=(a_0,a_1,a_2,\dots)\in \ell_p$ one has
$$\langle\phi_0^{\ast},x\rangle \phi_0=(a_0,0,0,\dots).$$
Let $(\{\psi_n\}_{n=0}^{\infty},\{\psi_n^{\ast}\}_{n=0}^{\infty})$ be a biorthogonal system in $\ell_p$, i.e. $\langle\psi_i^{\ast},\psi_j\rangle =\delta_{i}^{j}$ for $i,j=0,1,2,\dots$. Then the following statements hold

(i) If $p\in [1,2]$ and for every $x\in \ell_p$ we have
$$\left(\sum\limits_{n=1}^{\infty} |\langle\psi_n^{\ast}, x\rangle \langle\phi_n^{\ast}, \psi_n\rangle - \langle\phi_n^{\ast}, x\rangle  |^p \|\phi_n\|^p\right)^{\frac{1}{p}}\leq \varsigma_1 \|x\|,$$
where $\varsigma_1\in \left[0,(2M)^{-1}\right)$, then $\{\psi_n\}_{n=0}^{\infty}$ is also an unconditional basis of $\ell_p$, $p\in [1,2]$, isomorphic to $\{\phi_n\}_{n=0}^{\infty}$.

(ii)  If $p\in [2,\infty)$ and for every $x\in \ell_p$ we have
$$\left(\sum\limits_{n=1}^{\infty} |\langle\psi_n^{\ast}, x\rangle \langle\phi_n^{\ast}, \psi_n\rangle - \langle\phi_n^{\ast}, x\rangle  |^2 \|\phi_n\|^2\right)^{\frac{1}{2}}\leq \varsigma_2 \|x\|,$$
where $\varsigma_2\in \left[0,\frac{1}{\sqrt{8} M} \left(\frac{\Gamma \left( \frac{p+1}{2} \right)}{\sqrt{\pi}} \right)^{-\frac{1}{p}} \right)$, then $\{\psi_n\}_{n=0}^{\infty}$ is also an unconditional basis of $\ell_p$, $p\in [2,\infty)$, isomorphic to $\{\phi_n\}_{n=0}^{\infty}$.
\end{corollary}
Now we focus on the concept of symmetric basis of the space $E$ and give one sufficient condition for an arbitrary biorthogonal system in $\ell_p$ to be a symmetric basis of $\ell_p$.
\begin{definition}[\cite{Singer1,Singer3}]\label{symmetric basis}
A basis $\{\phi_n\}_{n=0}^{\infty}$ of $E$ is called symmetric provided
$$\sup_{\sigma \in \Pi} \sup_{|\beta_i| \leq 1, n\geq 0} \left\| \sum\limits_{i=0}^{n} \beta_i \langle \phi_i^{\ast},x\rangle \phi_{\sigma(i)} \right\| < \infty$$
holds for each $x\in E$, where $\{\phi_n^{\ast}\}_{n=0}^{\infty} \subset E^{\ast}$ is a sequence of coordinate functionals associated to $\{\phi_n\}_{n=0}^{\infty}$, and where $\Pi$ denotes the set of all permutations of the set $\mathbb{Z}_{+}.$
\end{definition}
The concept of symmetric basis was first introduced by Singer \cite{Singer3}.
For instance, vectors $\{e_n\}_{n=0}^{\infty}$ form a symmetric basis of the spaces $\ell_p$, $p\geq 1$, and $c_0$. Obviously, every symmetric basis is unconditional, but the converse is not true. Indeed, the normalized Haar basis in $L_p[0,1]$, $2\neq p>1$, is unconditional and, however, is not symmetric \cite{Singer1,Singer3}. It must be emphasized that the spaces $L_p[0,1]$, $2\neq p>1$, have no symmetric basis at all \cite{Singer1}. Concerning the question of uniqueness, up to isomorphism, of a symmetric basis, it is known that the spaces $\ell_p$, $p\geq 1$, and $c_0$ have a unique symmetric basis \cite{Lindenstrauss}. Combining this fact with Corollary \ref{ell p bases}, we obtain the following result.
\begin{proposition}\label{ell p symmetric bases}
Let $(\{\psi_n\}_{n=0}^{\infty},\{\psi_n^{\ast}\}_{n=0}^{\infty})$ be a biorthogonal system in $\ell_p$, i.e. $\langle\psi_i^{\ast},\psi_j\rangle =\delta_{i}^{j}$ for $i,j=0,1,2,\dots$. Then the following statements hold

(i) If $p\in [1,2]$ and for every $x\in \ell_p$ we have
$$\left(\sum\limits_{n=1}^{\infty} |\langle\psi_n^{\ast}, x\rangle \langle e_n^{\ast}, \psi_n\rangle - \langle e_n^{\ast}, x\rangle  |^p \right)^{\frac{1}{p}}\leq \varsigma_1 \|x\|,$$
where $\varsigma_1\in \left[0, \frac{1}{2}\right)$, then $\{\psi_n\}_{n=0}^{\infty}$ is a symmetric basis of $\ell_p$, $p\in [1,2]$.

(ii)  If $p\in [2,\infty)$ and for every $x\in \ell_p$ we have
$$\left(\sum\limits_{n=1}^{\infty} |\langle\psi_n^{\ast}, x\rangle \langle e_n^{\ast}, \psi_n\rangle - \langle e_n^{\ast}, x\rangle  |^2 \right)^{\frac{1}{2}}\leq \varsigma_2 \|x\|,$$
where $\varsigma_2\in \left[0,\frac{1}{\sqrt{8}} \left(\frac{\Gamma \left( \frac{p+1}{2} \right)}{\sqrt{\pi}} \right)^{-\frac{1}{p}} \right)$, then $\{\psi_n\}_{n=0}^{\infty}$ is a symmetric basis of $\ell_p$, $p\in [2,\infty)$.
\end{proposition}
\subsection{Isomorphic Schauder decompositions in spaces with\newline Schauder-Orlicz decompositions}
Let $E$ be a Banach space possessing a Schauder decomposition. We introduce the following definition.
\begin{definition}\label{Schauder-Orlicz decomposition}
A Schauder decomposition $\{\mathfrak{M}_n\}_{n=0}^{\infty}$ of $E$ with the a.s.c.p. $\{P_n\}_{n=0}^{\infty}$ will be said to be a Schauder-Orlicz decomposition provided there exists an Orlicz function $\Phi$ such that
\begin{equation}\label{Schauder-Orlicz}
    \inf\left\{\rho >0:\quad  \sum\limits_{n=0}^{\infty} \Phi\left(\frac{\|P_n x\|}{\rho}\right)\leq 1\right\} = \|x\| \quad \text{for all $x\in E$}.
\end{equation}
\end{definition}
For example, an Orlicz spaces $\ell_{\Phi}$, where $\Phi$ satisfies the $\Delta_2$-condition at zero, possess a Schauder-Orlicz decomposition. To see this it is sufficient to consider a Schauder decomposition $\{\mathfrak{M}_n\}_{n=0}^{\infty}$ corresponding to a basis $\{e_n\}_{n=0}^{\infty}$, i.e. $\{\mathfrak{M}_n=Lin \{e_n\}\}_{n=0}^{\infty}$. Also, by the generalized Parseval identity, we have that a Schauder decomposition in a Hilbert space $H$ is orthogonal if and only if it is a Schauder-Orlicz decomposition with an Orlicz function $\Phi(t)=t^2$. Thus, a concept of Schauder-Orlicz decomposition may be considered as natural generalization of the concept of an orthogonal Schauder decomposition  to the case of Banach spaces. For $\ell_{\Psi}$-Hilbertian ($\infty$-Hilbertian) Schauder decompositions $\{\mathfrak{M}_n\}_{n=0}^{\infty}$ in spaces with Schauder-Orlicz decompositions, we have the following result.
\begin{theorem}\label{MainTheorem2}
Let $E$ has a Schauder-Orlicz decomposition $\{\mathfrak{N}_n\}_{n=0}^{\infty}$ with the a.s.c.p. $\{F_n\}_{n=0}^{\infty}$ such that $\dim F_0 <\infty$, and let $\{\mathfrak{M}_n\}_{n=0}^{\infty}$ be an $\ell_{\Psi}$-Hilbertian ($\infty$-Hilbertian) Schauder decomposition of $E$ with constant $C$ and the a.s.c.p. $\{P_n\}_{n=0}^{\infty}$, where $P_0=F_0$. Assume that $\{J_n\}_{n=0}^{\infty}$ is a sequence of nonzero projections in $E$ satisfying $J_n J_m=\delta_{n}^{m} J_n$ for $n,m=0,1,\dots$ Furthermore, suppose that the condition (\ref{1}) holds and for all $x\in E$ we have
\begin{equation}\label{**}
    \left\|\sum\limits_{n=1}^{\infty} \|P_n (J_n-P_n) x\| e_n \right\|_{\ell_{\Psi}(c_0)} \leq \varsigma \|x\|,
\end{equation}
where $\varsigma$ is a constant such that $0\leq \varsigma<C^{-1}$.

Then $\{\mathfrak{M}_n\}_{n=0}^{\infty}$ and $\{J_n E\}_{n=0}^{\infty}$ are isomorphic Schauder decompositions of the space $E$.
\end{theorem}
We observe that, for the case of separable Orlicz spaces $E=\ell_{\Phi}$, Theorem \ref{MainTheorem1} follows from a Theorem \ref{MainTheorem2}. The condition on Schauder decomposition in $E$ to be an $\ell_{\Psi}$-Hilbertian, which is significant for Theorem \ref{MainTheorem1} and Theorem \ref{MainTheorem2}, is not very restrictive. Indeed, every Schauder decomposition is $1$-Hilbertian. Furthermore, every Schauder decomposition of uniformly convex Banach space $E$, which corresponds to a bounded basis $\{\phi_n\}_{n=0}^{\infty}$ of $E$, is $p$-Hilbertian for some $p>1$. More precisely, let $\{\phi_n\}_{n=0}^{\infty}$ be a bounded basis of uniformly convex Banach space $E$. Define a sequence of projections $\{P_n\}_{n=0}^{\infty}$ by $P_n x=\langle \phi_n^{\ast}, x \rangle \phi_n$, $x\in E$, $n=0,1,\dots,$ where $\{\phi_n^{\ast}\}_{n=0}^{\infty}$ is a sequence of coordinate functionals associated to $\{\phi_n\}_{n=0}^{\infty}$. Then there exist $p>1$, which depends on $\{\phi_n\}_{n=0}^{\infty}$ and on the modulus of convexity of a Banach space $E$, such that $\{P_n E\}_{n=0}^{\infty}$ is $p$-Hilbertian Schauder decomposition of $E$, see \cite{Gurarii} for details.

Also we note that the condition (\ref{1}) is an essential for Theorem \ref{MainTheorem1} and Theorem \ref{MainTheorem2}, as well as for Theorem \ref{Kato}, but, at the same time, is the only restriction imposed on the dimensions of projections $P_n$, $J_n$, $n\in Z_{+}$.
Further, since every Schauder decomposition is $1$-Hilbertian with constant $1$, we obtain the following consequence of Theorem \ref{MainTheorem2}, which is similar to the Corollary \ref{1-Hilbertian}, concerning stability of Schauder decompositions of certain structure in spaces possessing Schauder-Orlicz decompositions.
\begin{corollary}\label{1-Hilbertian-SO}
Let $E$ has a Schauder-Orlicz decomposition $\{\mathfrak{N}_n\}_{n=0}^{\infty}$ with the the a.s.c.p. $\{F_n\}_{n=0}^{\infty}$ such that $\dim F_0 <\infty$, and let  $\{\mathfrak{M}_n\}_{n=0}^{\infty}$ be a Schauder decomposition of $E$ with the a.s.c.p. $\{P_n\}_{n=0}^{\infty}$, where $P_0=F_0$. Assume that $\{J_n\}_{n=0}^{\infty}$ is a sequence of nonzero projections in $E$ such that $J_n J_m=\delta_{n}^{m} J_n$ for $n,m=0,1,\dots$ Furthermore, suppose that the condition (\ref{1}) holds and for all $x\in E$ we have
$$\sum\limits_{n=1}^{\infty} \|P_n (J_n-P_n) x\|\leq c \|x\|,$$
where $c\in [0,1)$. Then $\{\mathfrak{M}_n\}_{n=0}^{\infty}$ and $\{J_n E\}_{n=0}^{\infty}$ are isomorphic Schauder decompositions of $E$.
\end{corollary}
In case when all the subspaces of $\{\mathfrak{M}_n\}_{n=0}^{\infty}$ are one dimensional, by Corollary \ref{1-Hilbertian-SO}, we have the following result on stability of bases of certain structure in spaces possessing Schauder-Orlicz decompositions.
\begin{corollary}\label{1-Hilbertian bases}
Let $E$ has a Schauder-Orlicz decomposition $\{\mathfrak{M}_n\}_{n=0}^{\infty}$ with the a.s.c.p. $\{P_n\}_{n=0}^{\infty}$, where $\dim P_0=1$, and let
$\{\phi_n\}_{n=0}^{\infty}$ be a basis of $E$ with corresponding sequence of coordinate functionals $\{\phi_n^{\ast}\}_{n=0}^{\infty}$, such that for all $x\in E$ one has
$$\langle\phi_0^{\ast},x\rangle \phi_0=P_0 x.$$
Let $(\{\psi_n\}_{n=0}^{\infty},\{\psi_n^{\ast}\}_{n=0}^{\infty})$ be a biorthogonal system in $E$, i.e. $\langle\psi_i^{\ast},\psi_j\rangle =\delta_{i}^{j}$ for $i,j=0,1,2,\dots$. If for all $x\in E$ we have
$$\sum\limits_{n=1}^{\infty} |\langle\psi_n^{\ast}, x\rangle \langle\phi_n^{\ast}, \psi_n\rangle - \langle\phi_n^{\ast}, x\rangle  | \|\phi_n\|\leq c \|x\|,$$
where $c$ is a constant such that $0\leq c<1$, then $\{\phi_n\}_{n=0}^{\infty}$ and $\{\psi_n\}_{n=0}^{\infty}$ are isomorphic bases of the space $E$.
\end{corollary}
To show what does Theorem \ref{MainTheorem2} say for the case of concrete Banach spaces possessing Schauder-Orlicz decompositions we consider, for example, Schauder decompositions of certain structure in Hilbert spaces $H$ and unconditional Schauder decompositions in $H$. By virtue of Theorem \ref{MainTheorem2}, we have the following.
\begin{proposition}\label{Hilbert spaces}
Let $H$ be a Hilbert space with an orthogonal Schauder decomposition $\{\mathfrak{N}_n\}_{n=0}^{\infty}$ and the a.s.c.o. $\{F_n\}_{n=0}^{\infty}$ such that $\dim F_0 <\infty$, and let $\{\mathfrak{M}_n\}_{n=0}^{\infty}$ be a Schauder decomposition of $H$ with the a.s.c.p. $\{P_n\}_{n=0}^{\infty}$, where $P_0=F_0$. Assume that $\{J_n\}_{n=0}^{\infty}$ is a sequence of nonzero projections in $H$ such that $J_n J_m=\delta_{n}^{m} J_n$ for $n,m=0,1,\dots$ Furthermore, suppose that the condition (\ref{1}) holds and for all $x\in H$ we have
$$\sum\limits_{n=1}^{\infty} \|P_n (J_n-P_n) x\|\leq c \|x\|,$$
where $c$ is a constant such that $0\leq c<1$. Then $\{\mathfrak{M}_n\}_{n=0}^{\infty}$ and $\{J_n H\}_{n=0}^{\infty}$ are isomorphic Schauder decompositions of a Hilbert space $H$.
\end{proposition}
In case when all the subspaces of $\{\mathfrak{M}_n\}_{n=0}^{\infty}$ are one dimensional, Proposition \ref{Hilbert spaces} leads to the following.
\begin{corollary}\label{Hilbert space bases}
Let $\{h_n\}_{n=0}^{\infty}$ be an orthonormal basis of a Hilbert space $H$, and let $\{\phi_n\}_{n=0}^{\infty}$ be a basis of $H$ with corresponding sequence of coordinate functionals $\{\phi_n^{\ast}\}_{n=0}^{\infty}$, such that $\phi_0=\phi_0^{\ast}=h_0$.
Let $(\{\psi_n\}_{n=0}^{\infty},\{\psi_n^{\ast}\}_{n=0}^{\infty})$ be a biorthogonal system in $H$, i.e. $\langle\psi_i^{\ast},\psi_j\rangle =\delta_{i}^{j}$ for $i,j=0,1,2,\dots$. If for all $x\in H$ we have
$$\sum\limits_{n=1}^{\infty} |\langle\psi_n^{\ast}, x\rangle \langle\phi_n^{\ast}, \psi_n\rangle - \langle\phi_n^{\ast}, x\rangle  | \|\phi_n\|\leq c \|x\|,$$
where $c$ is a constant such that $0\leq c<1$, then $\{\phi_n\}_{n=0}^{\infty}$ and $\{\psi_n\}_{n=0}^{\infty}$ are isomorphic bases of a Hilbert space $H$.
\end{corollary}
Applying Theorem \ref{MainTheorem2}, in view of Remark \ref{p=q=2}, we obtain the following results on stability of unconditional Schauder decompositions and unconditional bases in $H$.
\begin{theorem}\label{Hilbert spaces unconditional}
Let $H$ be a Hilbert space with an orthogonal Schauder decomposition $\{\mathfrak{N}_n\}_{n=0}^{\infty}$ and the a.s.c.o. $\{F_n\}_{n=0}^{\infty}$ such that $\dim F_0 <\infty$, and let $\{\mathfrak{M}_n\}_{n=0}^{\infty}$ be an unconditional Schauder decomposition of $H$ with constant $M$ and the a.s.c.p. $\{P_n\}_{n=0}^{\infty}$, where $P_0=F_0$. Assume that $\{J_n\}_{n=0}^{\infty}$ is a sequence of nonzero projections in $H$ such that $J_n J_m=\delta_{n}^{m} J_n$ for $n,m=0,1,\dots$
Furthermore, suppose that the condition (\ref{1}) holds and for all $x\in H$ we have
$$\left(\sum\limits_{n=1}^{\infty} \|P_n (J_n-P_n) x\|^2 \right)^{\frac{1}{2}}\leq \varsigma \|x\|,$$
where $\varsigma$ is a constant such that $0\leq \varsigma<(2M)^{-1}$.

Then $\{J_n H\}_{n=0}^{\infty}$ is also an unconditional Schauder decomposition of $H$, isomorphic to $\{\mathfrak{M}_n\}_{n=0}^{\infty}$.
\end{theorem}
\begin{corollary}\label{Hilbert space unconditional bases}
Let $\{h_n\}_{n=0}^{\infty}$ be an orthonormal basis of a Hilbert space $H$, and let $\{\phi_n\}_{n=0}^{\infty}$ be an unconditional basis of $H$ with constant $M$ and corresponding sequence of coordinate functionals $\{\phi_n^{\ast}\}_{n=0}^{\infty}$, such that $\phi_0=\phi_0^{\ast}=h_0$.
Let $(\{\psi_n\}_{n=0}^{\infty},\{\psi_n^{\ast}\}_{n=0}^{\infty})$ be a biorthogonal system in $H$, i.e. $\langle\psi_i^{\ast},\psi_j\rangle =\delta_{i}^{j}$ for $i,j=0,1,2,\dots$. If for all $x\in H$ we have
$$\left(\sum\limits_{n=1}^{\infty} |\langle\psi_n^{\ast}, x\rangle \langle\phi_n^{\ast}, \psi_n\rangle - \langle\phi_n^{\ast}, x\rangle  |^2 \|\phi_n\|^2 \right)^{\frac{1}{2}}\leq \varsigma \|x\|,$$
where $\varsigma\in [0,(2M)^{-1})$, then $\{\psi_n\}_{n=0}^{\infty}$ is also an unconditional basis of $H$, isomorphic to $\{\phi_n\}_{n=0}^{\infty}$.
\end{corollary}
Since by Proposition \ref{ell_1,c_0}, $\textit{(ii)}$, every unconditional Schauder decomposition of the space $c_0$ is a $c_0$-decomposition, for each unconditional Schauder decomposition $\{\mathfrak{M}_n\}_{n=0}^{\infty}$ of the space $c_0$ with the a.s.c.p. $\{P_n\}_{n=0}^{\infty}$ there exists a constant $C>0$ such that
\begin{equation}\label{c_0}
    \left\|x\right\| \leq C \max\limits_{n\geq 0} \|P_n x\|\quad \text{for every $x\in c_0$}.
\end{equation}
Therefore, Theorem \ref{MainTheorem2} yields the following stability result for unconditional Schauder decompositions in $c_0$.
\begin{proposition}\label{c_0 unconditional}
Let $\{\mathfrak{N}_n\}_{n=0}^{\infty}$ be a Schauder-Orlicz decomposition of $c_0$ with the a.s.c.p. $\{F_n\}_{n=0}^{\infty}$ such that $\dim F_0 <\infty$, and
let $\{\mathfrak{M}_n\}_{n=0}^{\infty}$ be an unconditional Schauder decomposition of the space $c_0$ with the a.s.c.p. $\{P_n\}_{n=0}^{\infty}$ satisfying (\ref{c_0}), where $P_0=F_0$.
Assume that $\{J_n\}_{n=0}^{\infty}$ is a sequence of nonzero projections in $c_0$ such that $J_n J_m=\delta_{n}^{m} J_n$ for $n,m=0,1,\dots$
Furthermore, suppose that the condition (\ref{1}) holds and for all $x\in c_0$ we have
$$\max\limits_{n\geq 0} \|P_n (J_n-P_n) x\|\leq \varsigma \|x\|,$$
where $\varsigma$ is a constant such that $0\leq \varsigma<C^{-1}$.

Then $\{J_n (c_0)\}_{n=0}^{\infty}$ is also an unconditional Schauder decomposition of the space $c_0$, isomorphic to $\{\mathfrak{M}_n\}_{n=0}^{\infty}$.
\end{proposition}

\section{Proofs}
\subsection{Proof of Proposition \ref{proposition}}
Let (i) holds, i.e. there exists an isomorphism $S$ such that $\mathfrak{N}_n=S\mathfrak{M}_n$, $n=0,1,\dots$ Then each fixed $x\in E$ has an expansion $ x=\sum\limits_{n=0}^{\infty} J_n x = \sum\limits_{n=0}^{\infty} x_n$, $x_n\in \mathfrak{N}_n,$ $n=0,1,\dots$ Thus, an element $y=S^{-1}x\in E$ has an expansion
\begin{equation}\label{prop1}
    y=\sum\limits_{n=0}^{\infty} S^{-1} x_n,\: x_n\in \mathfrak{N}_n,\: n=0,1,\dots
\end{equation}
On the other hand we observe that $y$ has also the following expansion
\begin{equation}\label{prop2}
    y=\sum\limits_{n=0}^{\infty} y_n=\sum\limits_{n=0}^{\infty} P_n y,\: y_n\in \mathfrak{M}_n,\: n=0,1,\dots
\end{equation}
Since $\mathfrak{M}_n=S^{-1}\mathfrak{N}_n$ for all $n$, we have that $S^{-1} x_n \in \mathfrak{M}_n$, $n=0,1,\dots$ Combining (\ref{prop1}) with (\ref{prop2}) we conclude that $S^{-1} x_n=P_n y$ and $x_n=J_n x=S P_n S^{-1}x,$ $n=0,1,\dots$

Conversely, if (ii) satisfied, i.e. if there exists an isomorphism $S$ such that $J_n=S P_n S^{-1}$ for each $n$, then, for every fixed $n$ and $x\in \mathfrak{N}_n$, $ x=J_nx=SP_nS^{-1} x=S y$, where $y=P_nS^{-1} x \in \mathfrak{M}_n$. Hence, $x\in S \mathfrak{M}_n$ and $\mathfrak{N}_n \subset S \mathfrak{M}_n.$ To prove the converse inclusion, consider any fixed $n$ and $y\in \mathfrak{M}_n$. Then $P_n y=y$ and $x=Sy=SP_n y\in S \mathfrak{M}_n$. Since $J_n x=J_n S P_n y=SP_n y=x$ we obtain that $x\in \mathfrak{N}_n$.
\subsection{Proof of Lemma \ref{Lemma}}
Since $E$ has or Orlicz-Rademacher type $\Phi$, or Orlicz-Rademacher infratype $\Phi$, for any element $x\in E$ and for every finite set of vectors $\{P_j x\}_{j=0}^{n}\subset E$ there exists a set of numbers $\{\underline{\varepsilon}_j \}_{j=0}^{n}\subset \{-1,1\}$ such that
$$\left\| \sum\limits_{j=0}^{n} \underline{\varepsilon}_j P_j x \right\| = \left(\min_{\varepsilon_j=\pm 1} \left\| \sum\limits_{j=0}^{n} \varepsilon_j P_j x \right\|^2\right)^{\frac{1}{2}} $$
$$\leq\left \{
\begin{array}{l}
\displaystyle \left(\mathbb{E} \left\| \sum\limits_{j=0}^{n} \varepsilon_j P_j x \right\|^2\right)^{\frac{1}{2}} \leq T_{\Phi}(E) \inf\left\{\rho >0:\quad  \sum\limits_{j=0}^{n} \Phi\left(\frac{\|P_j x\|}{\rho}\right)\leq 1\right\}\\[10pt]\displaystyle
I_{\Phi}(E) \inf\left\{\rho >0:\quad  \sum\limits_{j=0}^{n} \Phi\left(\frac{\|P_j x\|}{\rho}\right)\leq 1\right\}\\
\end{array}
\right .$$
$\leq T(\Phi) \inf\left\{\rho >0:\quad  \sum\limits_{j=0}^{n} \Phi\left(\frac{\|P_j x\|}{\rho}\right)\leq 1\right\}$. Construct the operators $P_n^{+}=\sum\limits_{j: \underline{\varepsilon}_j=1} P_j$, $P_n^{-}=\sum\limits_{j: \underline{\varepsilon}_j=-1} P_j$. Then,
$$\|x\|= \lim_{n\rightarrow \infty} \left\| (P_n^{+}-P_n^{-})^2 x \right\| \leq 2 M \lim_{n\rightarrow \infty} \left\| (P_n^{+}-P_n^{-}) x \right\| $$
$$\leq 2 M T(\Phi)\inf\left\{\rho >0:\quad  \sum\limits_{j=0}^{\infty} \Phi\left(\frac{\|P_j x\|}{\rho}\right)\leq 1\right\}$$
and the right inequality from (\ref{Type-Cotype}) is proved with $T=2M T(\Phi)$.

Further, since $E$ has or Orlicz-Rademacher cotype $\Psi$, or Orlicz-Rademacher $M$-cotype $\Psi$, for any element $x\in E$ and for every finite set of vectors $\{P_j x\}_{j=0}^{n}\subset E$ there exists a set of numbers $\{\overline{\varepsilon}_j \}_{j=0}^{n}\subset \{-1,1\}$ such that $$\left\| \sum\limits_{j=0}^{n} \overline{\varepsilon}_j P_j x \right\| = \left(\max_{\varepsilon_j=\pm 1} \left\| \sum\limits_{j=0}^{n} \varepsilon_j P_j x \right\|^2\right)^{\frac{1}{2}} $$
$$\geq\left \{
\begin{array}{l}
\displaystyle \left(\mathbb{E} \left\| \sum\limits_{j=0}^{n} \varepsilon_j P_j x \right\|^2\right)^{\frac{1}{2}} \geq C_{\Psi}(E) \inf\left\{\rho >0:\quad  \sum\limits_{j=0}^{n} \Psi\left(\frac{\|P_j x\|}{\rho}\right)\leq 1\right\}\\[10pt]\displaystyle
M_{\Psi}(E) \inf\left\{\rho >0:\quad  \sum\limits_{j=0}^{n} \Psi\left(\frac{\|P_j x\|}{\rho}\right)\leq 1\right\}\\
\end{array}
\right .$$
$\geq C(\Psi) \inf\left\{\rho >0:\quad  \sum\limits_{j=0}^{n} \Psi\left(\frac{\|P_j x\|}{\rho}\right)\leq 1\right\}$. We observe that, for each set of numbers $\{\overline{\varepsilon}_j \}_{j=0}^{n}\subset \{-1,1\}$ there exist two sets of numbers $\{\delta_j^{+} \}_{j=0}^{n}\subset \{0,1\}$ and $\{\delta_j^{-} \}_{j=0}^{n}\subset \{0,1\}$ such that
$$\left\| \sum\limits_{j=0}^{n} \overline{\varepsilon}_j P_j x \right\|=\left\| \sum\limits_{j=0}^{n} \delta_j^{+} P_j x - \sum\limits_{j=0}^{n} \delta_j^{-} P_j x \right\|\leq 2M \|x\|.$$
Consequently, $C(\Psi) \inf\left\{\rho >0:\quad  \sum\limits_{j=0}^{n} \Psi\left(\frac{\|P_j x\|}{\rho}\right)\leq 1\right\}\leq 2M \|x\|$, and (\ref{Type-Cotype}) is proved with $C=(2M)^{-1} C(\Psi)$.
\subsection{Proof of Theorem \ref{MainTheorem1}}
Construct the operator $S$ on the Orlicz space $\ell_{\Phi}$ by
\begin{equation}\label{S}
    S=\sum\limits_{n=0}^{\infty} P_n J_n.
\end{equation}
To show that $S$ exists in the strong sense we show that
$$\sum\limits_{n=0}^{\infty} \left( P_n -P_n J_n \right)=\sum\limits_{n=0}^{\infty} P_n  \left( P_n -J_n \right)$$
converges strongly. Indeed, since $\{\mathfrak{M}_n\}_{n=0}^{\infty}$ is an $\ell_{\Psi}$-Hilbertian ($\infty$-Hilbertian) Schauder decomposition of $\ell_{\Phi}$ with constant $C$ and the a.s.c.p. $\{P_n\}_{n=0}^{\infty}$, for every $x\in\ell_{\Phi}$ and for each $N\in\mathbb{Z}_{+}$ we have by (\ref{**}) that
$$\left\|\sum\limits_{n=k}^{k+N} P_n  \left( P_n -J_n \right) x \right\| \leq C \left\|\sum\limits_{j=0}^{\infty} \left\|P_j \left(\sum\limits_{n=k}^{k+N} P_n  \left( P_n -J_n \right) x\right)\right\| e_j\right\|_{\ell_{\Psi}(c_0)}$$
$$= C \left\|\sum\limits_{n=k}^{k+N} \|P_n (J_n-P_n) x\| e_n \right\|_{\ell_{\Psi}(c_0)}$$
tends to zero when $k$ tends to infinity. Therefore, $\sum\limits_{n=0}^{\infty} P_n  \left( P_n -J_n \right) x$ is convergent and, consequently, the series
$$\sum\limits_{n=0}^{\infty} P_n J_n x = \sum\limits_{n=0}^{\infty} P_n x - \sum\limits_{n=0}^{\infty} P_n \left( P_n -J_n \right) x$$
is also convergent. Consider the operator $R= \sum\limits_{n=1}^{\infty} P_n  \left( P_n -J_n \right) = I - P_0^{} - \sum\limits_{n=1}^{\infty} P_n J_n$ and note that, since
$$ \|Rx\|= \left\|\sum\limits_{n=1}^{\infty} P_n  \left( P_n -J_n \right) x \right\| \leq C \left\|\sum\limits_{j=0}^{\infty} \left\|P_j \left(\sum\limits_{n=1}^{\infty} P_n  \left( P_n -J_n \right) x\right)\right\| e_j\right\|_{\ell_{\Psi}(c_0)}$$
$$\leq C\varsigma \|x\|$$
by (\ref{**}), we have that $\|R\|<1.$ Further we observe that, since $S=P_0^{} J_0^{} +I -P_0^{} - R$, $\|S\| < \|J_0^{}\|+3<\infty$. Now (\ref{S}) implies that
$$
S J_n=\sum\limits_{j=0}^{\infty} P_j J_j  J_n =  P_n J_n =  \sum\limits_{j=0}^{\infty} P_n P_j J_j =P_n S.$$
Thus the theorem will be proved if we show that $S$ is continuously invertible. To this end we consider
\begin{equation}\label{tildeS}
    \widetilde{S}=\sum\limits_{n=1}^{\infty} P_n J_n = I -P_0^{} -R.
\end{equation}
Since  $\dim P_0^{}=m <\infty$, by the definition of projection $P_0$ we have that $\left( I-P_0^{}\right)$  is a Fredholm operator with
$$nul \left( I-P_0^{}\right) = m,\quad ind \left( I-P_0^{}\right) =0,\quad \gamma \left( I-P_0^{}\right)=1,$$
where $nul \: T$ denotes the nullity,  $ind \: T$ the index, and $\gamma\: (T)$ the reduced\newline minimum modulus, of the operator $T$ (for these notions see, e.g., \cite{Kato1}, Chapter IV, $\S 5.1$). Indeed, first we note that $nul \left( I-P_0^{}\right) = \dim P_0^{} =m,$
$$def \left( I-P_0^{}\right) = \dim \ell_{\Phi}|_{Im \left( I-P_0^{}\right) }= \dim \ell_{\Phi}|_{\overline{{Im \left( I-P_0^{}\right) }}}= \dim co\ker\left( I-P_0^{}\right)$$ 
$$= \dim \left( Im \left( I-P_0^{}\right)  \right)^{\bot}=m,$$
$ind \left( I-P_0^{}\right)= nul \left( I-P_0^{}\right) - def \left( I-P_0^{}\right) =0,$
where $def \: T$ denotes the deficiency of $T$, see, e.g, \cite{Bilalov2,Kato1}. Second, since for each $x=(a_0,a_1,a_2,\dots)\in \ell_{\Phi}$,
$$\inf\limits_{v\in \ker \left( I-P_0^{}\right)} \|x-v\|=\inf\limits_{v\in Im P_0^{}} \inf\left\{\rho >0:\quad  \sum\limits_{n=0}^{\infty} \Phi\left(\frac{|a_n-v_n|}{\rho}\right)\leq 1\right\}$$
$$=\inf\left\{\rho >0:\quad  \sum\limits_{n=m}^{\infty} \Phi\left(\frac{|a_n|}{\rho}\right)\leq 1\right\}$$
$=\|\left( I-P_0^{}\right) x\|$, where $v=(v_0,v_1,v_2,\dots)\in Im P_0$, we obtain that
$$\gamma \left( I-P_0^{}\right)=$$ 
$$\sup \left\{ \gamma: \: \|\left( I-P_0^{}\right) x\| \geq \gamma \inf\limits_{v\in \ker \left( I-P_0^{}\right)} \|x-v\|, \: x\in D\left( I-P_0^{}\right)= \ell_{\Phi} \right\}=1.$$

Since $\|R\|<1= \gamma \left( I-P_0^{}\right)$, it follows that $\widetilde{S}=\left(I-P_0^{}\right)-R$ is also Fredholm, with
\begin{equation}\label{chtildeS}
    nul\: \widetilde{S} \leq nul \left( I-P_0^{}\right) =m, \quad\quad\quad ind\: \widetilde{S} = ind \left( I-P_0^{}\right)=0
\end{equation}
(see \cite{Kato1}, Chapter IV, Theorem 5.22). Since $S=P_0^{} J_0^{}+ \widetilde{S}$, where $P_0^{} J_0^{}$ is compact, $S$ is also Fredholm and $ind \: S=ind \: \widetilde{S}=0$ (see \cite{Kato1}, Chapter IV, Theorem 5.26). Consequently, $nul\: S =def\: S$, and $S$ will be invertible on $\ell_{\Phi}$ if and only if $nul\: S =def\: S=0$. Therefore it is sufficient to show that $nul\: S=0$.

To this end we first show that
\begin{equation}\label{Ker}
    \ker \widetilde{S}= Im J_0.
\end{equation}
If $x\in Im J_0^{}$, i.e. $x=J_0^{} y$, then $\widetilde{S} x =\widetilde{S} J_0^{} y =  \sum\limits_{n=1}^{\infty} P_n J_n J_0^{} y=0$ and, consequently, $x\in \ker \widetilde{S}$. On the other hand, $\ker \widetilde{S} \subset Im J_0$, because $\ker \widetilde{S}$ and $Im J_0$ are linear subspaces, $\dim Im J_0^{}=m$ and $\dim \ker \widetilde{S}\leq m$ by (\ref{chtildeS}).

Assume now that $x\in \ker S$. Then,
$$0=P_0^{} S x = P_0^{} \sum\limits_{n=0}^{\infty} P_n J_n x = P_0 J_0 x $$
and $\widetilde{S}x =Sx-P_0 J_0 x =0$. Consequently,  $x\in \ker \widetilde{S}$, $x=J_0^{} y$ by (\ref{Ker}) and, hence, $$P_0 x=P_0 J_0 y=P_0^{} \sum\limits_{n=0}^{\infty} P_n J_n J_0^{} y =P_0^{} \sum\limits_{n=0}^{\infty} P_n J_n x=0.$$ As a result,
$$(I-R) x = (\widetilde{S}+ P_0^{})x =0.$$
Since $\|R\|<1$, we obtain $x=0$. Thus, $\ker S =\{0\}$, $nul\:S=0$ and $S$ is continuously invertible. To complete the proof we apply Proposition \ref{proposition}.
\subsection{Proof of Theorem \ref{MainTheorem2}}
The proof of this theorem is similar to the proof of Theorem \ref{MainTheorem1}. We only note that, since $E$ has a Schauder-Orlicz decomposition $\{\mathfrak{N}_n\}_{n=0}^{\infty}$ with the the a.s.c.p. $\{F_n\}_{n=0}^{\infty}$, then there exists an Orlicz function $\Phi$ such that for each $x\in E$ one has
$$\inf\limits_{v\in \ker \left( I-P_0^{}\right)} \|x-v\|=\inf\limits_{v\in Im F_0^{}} \inf\left\{\rho >0:\quad  \sum\limits_{n=0}^{\infty} \Phi\left(\frac{\|F_n(x-v)\|}{\rho}\right)\leq 1\right\}  $$
$$=\inf\left\{\rho >0:\quad  \sum\limits_{n=1}^{\infty} \Phi\left(\frac{\|F_n(x-F_0 x)\|}{\rho}\right)\leq 1\right\}$$
$$=\inf\left\{\rho >0:\quad  \sum\limits_{n=0}^{\infty} \Phi\left(\frac{\|F_n(x-F_0 x)\|}{\rho}\right)\leq 1\right\}=\|(I-P_0) x\|.$$
Therefore, as in the proof of Theorem \ref{MainTheorem1}, $$\gamma \left( I-P_0^{}\right) =$$ 
$$\sup \left\{ \gamma: \: \|\left( I-P_0^{}\right) x\| \geq \gamma \inf\limits_{v\in \ker \left( I-P_0^{}\right)} \|x-v\|, \: x\in D\left( I-P_0^{}\right)= E \right\}=1.$$
\section{Conclusions}
We obtain some stability theorems for Schauder decompositions in certain Banach spaces. These theorems may be considered as an extensions of Theorem \ref{Kato} of Kato on similarity for sequences of projections in Hilbert spaces to the case of Banach spaces, and may be classified as stability theorems in terms of the closeness of projections. More precisely, we introduce the class of Schauder-Orlicz decompositions and find conditions providing an existence of an isomorphism between certain Schauder decomposition and some sequence of nonzero subspaces. Stability theorems were obtained in the case of
Orlicz sequence spaces with Schauder decompositions (Theorem \ref{MainTheorem1}) and in the case of spaces possessing Schauder-Orlicz decompositions (Theorem \ref{MainTheorem2}).

We introduce the generalized concepts of type, cotype, infratype, $M$-cotype of a Banach space $E$ (Definitions \ref{O-type}, \ref{O-cotype}, \ref{O-infratype}, \ref{O-M-cotype}) and find interconnections between these geometric characteristics of $E$ and the properties of unconditional Schauder decompositions in $E$. As an application, we deduce one stability theorem of geometric type for unconditional Schauder decompositions in a Banach space (Theorem \ref{IsomorphicUSD}). The results are accompanied by a number of consequences concerning unconditional Schauder decompositions, isomorphic unconditional Schauder decompositions and isomorphic (equivalent) bases in spaces $L_p(\mu)$,  $\ell_p$ ($p\geq 1$), $c_0$ and in Hilbert spaces.

Just as Theorem \ref{Kato} of Kato plays a special role in the analysis of spectral properties of nonselfadjoint and unbounded linear operators in Hilbert spaces (see, e.g., \cite{Adduci1,Adduci2,Clark,Hughes,Kato1}), stability Theorems \ref{IsomorphicUSD}, \ref{MainTheorem1} and \ref{MainTheorem2} may be very useful in the study of spectral properties of unbounded linear operators acting on certain Banach spaces.

Concerning Schauder-Orlicz decompositions and Theorem \ref{MainTheorem2}, the following question becomes an important and requires further study. How wide the class of spaces possessing Schauder-Orlicz decompositions is? In particular: Does the space $L_p(\mu)$ possesses a Schauder-Orlicz decomposition  for some $2\neq p\geq1$? If the answer is positive, we can use Proposition \ref{Proposition} to obtain a stability theorem on unconditional Schauder decompositions in $L_p(\mu)$ spaces, similar to the Theorem \ref{ell_p Theorem}.

Concerning Proposition \ref{Proposition} it must be also emphasized that the constants in (\ref{1<p<p_0 inequality}), (\ref{p_0<p<2 inequality}), and (\ref{2<p<infty inequality}) are quite sharp. Nevertheless, the question on further improvement of the constants in (\ref{1<p<p_0 inequality}), (\ref{p_0<p<2 inequality}), and (\ref{2<p<infty inequality}) makes sense.


\begin{thebibliography}{9}

\bibitem{Adduci1} Adduci J., Mityagin B., Eigensystem of an $L^2$-perturbed harmonic oscillator is an unconditional basis, Cent. Eur. J. Math., 2012,
10(2), 569--589

\bibitem{Adduci2} Adduci J., Mityagin B., Root system of a perturbation of a selfadjoint operator with discrete spectrum, Integral Equations Operator
Theory, 2012, 73(2), 153--175

\bibitem{Ahmad} Ahmad K., A note on equivalence of sequences of subspaces in Banach spaces, An. Stiint. Univ. "Ovidius" Constanta Ser. Mat.,
1989, 27, 9--12

\bibitem{Allexandrov} Allexandrov G., Kutzarova D., Plichko A., A separable space with no Schauder decomposition, Proc. Amer. Math. Soc., 1999,
127(9), 2805--2806

\bibitem{Bilalov2} Bilalov B.T., Veliev S.G., Some Questions of Bases, Elm, Baku, 2010 (in Russian)

\bibitem{Bonet} Bonet J., Ricker W.J., Schauder decompositions and the Grothendieck and Dunford-Pettis properties in K\"{o}the echelon spaces of infinite order, Positivity, 2007, 11(1), 77--93

\bibitem{Casazza} Casazza P.G., Kalton N.J., Unconditional bases and unconditional finite-dimensional decompositions in \newline Banach spaces,
Israel J. Math., 1996, 95(1), 349--373

\bibitem{Chadwick} Chadwick J.J.M., Cross R.W., Schauder decompositions in non-separable Banach spaces, Bull. Aust. Math. Soc., 1972, 6(1), 133--144

\bibitem{Clark} Clark C., On relatively bounded perturbations of ordinary differential operators, Pacific J. Math., 1968, 25(1), 59--70

\bibitem{Clement} Clement P., De Pagter B., Sukochev F.A., Witvliet H., Schauder decompositions and multiplier theorems, Studia Math., 2000,
138(2), 135--163

\bibitem{Curtain} Curtain R.F., Zwart H.J., An Introduction to Infinite-Dimensional Linear Systems Theory, Texts in Applied Mathematics, Volume 21, Springer-Verlag, New-York, 1995

\bibitem{Davis} Davis W.J., Schauder decompositions in Banach spaces, Bull. Amer. Math. Soc., 1968, 74(6), 1083--1085

\bibitem{DelaRosa} De la Rosa M., Frerick L., Grivaux S., Peris A., Frequent hypercyclicity, chaos, and unconditional Schauder decompositions,
Israel J. Math., 2012, 190(1), 389--399

\bibitem{DePagter} De Pagter B., Ricker W.J., Products of commuting Boolean algebras of projections and Banach space \newline geometry,
Proc. Lond. Math. Soc. (3), 2005, 91(3), 483--508

\bibitem{Fage1} Fage M.K., Idempotent operators and their rectification, Dokl. Akad. Nauk, 1950, 73, 895--897 (in Russian)

\bibitem{Fage2} Fage M.K., The rectification of bases in Hilbert space, Dokl. Akad. Nauk, 1950, 74, 1053--1056 (in Russian)

\bibitem{Gohberg} Gohberg I.C., Krein M.G., Introduction to the Theory of Linear Nonselfadjoint Operators in Hilbert Space, Transl. Math. Monogr., 18, American Mathematical Society, Providence, Rhode Island, 1969

\bibitem{Grinblyum} Grinblyum M.M., On the representation of a space of type B in the form of a direct sum of subspaces, Dokl. Akad. Nauk, 1950,
 70, 749--752 (in Russian)

\bibitem{Gurarii} Gurarii V.I., Gurarii N.I., Bases in uniformly convex and uniformly smooth Banach spaces, Izv. Ross. Akad. Nauk Ser. Mat., 1971,
35(1), 210–-215 (in Russian)

\bibitem{Haagerup} Haagerup U., The best constants in the Khintchine inequality, Studia Math., 1982, 70, 231--283

\bibitem{Hughes} Hughes E., Perturbation theorems for relative spectral problems, Canad. J. Math., 1972, 24(1), 72--81

\bibitem{Jain1} Jain P.K,, Ahmad K., Maskey S.M., Domination and equivalence of sequences of subspaces in dual spaces, Czechoslovak Math. J.,
1986, 36(3), 351--357

\bibitem{Jain2} Jain P.K,, Ahmad K., Schauder decompositions and best approximations in Banach spaces, Port. Math., 1987, 44(1), 25--39

\bibitem{Jain3} Jain P.K,, Ahmad K., Unconditional Schauder decompositions and best approximations in Banach spaces, Indian J. Pure Appl. Math.,
1981, 12(12), 1456--1467

\bibitem{Johnson1} Johnson W.B., Finite-dimensional Schauder decompositions in $\pi_{\lambda}$ and dual $\pi_{\lambda}$ spaces, Illinois J. Math.,
1970, 14(4), 642--647

\bibitem{Johnson2} Johnson W.B., Lindenstrauss J., Handbook of the Geometry of Banach Spaces, Volume 1, Elsevier, 2001

\bibitem{Johnson3} Johnson W.B., Lindenstrauss J., Handbook of the Geometry of Banach Spaces, Volume 2, Elsevier, 2003

\bibitem{Kadets} Kadets M.I, Kadets V.M., Series in Banach Spaces, Conditional and Unconditional Convergence, Birkhauser, Berlin, 1997

\bibitem{Kato1} Kato T., Perturbation Theory for Linear Operators, 2nd ed. (reprint), Classics Math., Springer, Berlin, 1995

\bibitem{Kato2} Kato T., Similarity for sequences of projections, Bull. Amer. Math. Soc., 1967, 73(6), 904--905

\bibitem{Lindenstrauss} Lindenstrauss J., Tzafriri L., Classical Banach Spaces I and II, Reprint of the 1977, 1979 ed., Springer-Verlag, Berlin, 1996

\bibitem{Marcus} Marcus A.S., A basis of root vectors of a dissipative operator, Dokl. Akad. Nauk, 1960, 132(3), 524--527 (in Russian)

\bibitem{Rabah1} Rabah R., Sklyar G.M., Rezounenko A.V., Generalized Riesz basis property in the analysis of neutral type systems, C. R. Math. Acad. Sci. Paris, Ser. I, 2003, 337, 19--24

\bibitem{Rabah2} Rabah R., Sklyar G.M., Rezounenko A.V., Stability analysis of neutral type systems in Hilbert space, J. Differential Equations, 2005, 214, 391--428

\bibitem{Rabah3} Rabah R., Sklyar G.M., The analysis of exact controllability of neutral-type systems by the moment problem approach, SIAM J. Control Optim., 2007, 46(6), 2148--2181

\bibitem{Retherford1} Retherford J.R., Basic sequences and the Paley-Wiener criterion, Pacific J. Math., 1964, 14, 1019--1027

\bibitem{Retherford3} Retherford J.R., Some remarks on Schauder bases of subspaces, Rev. Roumaine Math. Pures Appl., 1966, 11, 787--792

\bibitem{Sanders1} Sanders B.L., Decompositions and reflexivity in Banach spaces, Proc. Amer. Math. Soc., 1965, 16(2), 204--208

\bibitem{Sanders2} Sanders B.L., On the existance of [Schauder] decompositions in Banach spaces, Proc. Amer. Math. Soc., 1965, 16(5), 987--990

\bibitem{Singer1} Singer I., Bases in Banach Spaces I, Springer-Verlag, Berlin, 1970

\bibitem{Singer2} Singer I., Bases in Banach Spaces II, Springer-Verlag, Berlin, 1981

\bibitem{Singer3} Singer I., On Banach spaces with symmetric bases, Rev. Roumaine Math. Pures Appl., 1961, 6, 159--166 (in Russian)

\bibitem{Wermer} Wermer J., Commuting spectral measures on Hilbert space, Pacific J. Math., 1954, 4, 355--361

\bibitem{Zwart} Zwart H., Riesz basis for strongly continuous groups, J. Differential Equations, 2010, 249, 2397--2408

\end{thebibliography}
\end{document}